\documentclass[10pt,a4paper]{article}

\usepackage{amsmath}
\usepackage{amsthm}
\usepackage{amssymb}
\usepackage{amscd}

\theoremstyle{definition}
\newtheorem{theorem}{Theorem}[section]

\newtheorem{lemma}[theorem]{Lemma}

\newtheorem{remark}[theorem]{Remark}
\newtheorem{conjecture}[theorem]{Conjecture}

\newenvironment{demo}[1]{%
  \trivlist
  \item[\hskip\labelsep
        {\bf #1.}]
}{%
\hfill\qedsymbol
  \endtrivlist
}

\newenvironment{roster}[1]{%
\begin{list}{}{%
\setlength{\topsep}{0pt}
\setlength{\itemsep}{0pt}
\setlength{\parsep}{0pt}
\settowidth{\labelwidth}{#1}
\setlength{\leftmargin}{\labelwidth}
\addtolength{\leftmargin}{\parindent}
}}{\end{list}}

\renewcommand\tilde{\widetilde}
\newcommand\Pf{\operatorname{Pf}}
\newcommand\perm{\operatorname{perm}}
\newcommand\Hf{\operatorname{Hf}}
\newcommand\rank{\operatorname{rank}}
\newcommand\sgn{\operatorname{sgn}}
\newcommand\Comp{\mathbb{C}}

\newcommand\Int{\mathbb{Z}}

\newcommand{\vectx}{\overrightarrow{\boldsymbol{x}}}
\newcommand{\vecty}{\overrightarrow{\boldsymbol{y}}}
\newcommand{\vectz}{\overrightarrow{\boldsymbol{z}}}

\newcommand{\vecta}{\overrightarrow{\boldsymbol{a}}}
\newcommand{\vectb}{\overrightarrow{\boldsymbol{b}}}
\newcommand{\vectc}{\overrightarrow{\boldsymbol{c}}}

\newcommand{\vectone}{\overrightarrow{\boldsymbol{1}}}
\newcommand{\vectzero}{\overrightarrow{\boldsymbol{0}}}
\newcommand\ASM{{\mathcal{A}}}
\newcommand\HTS{{\mathrm{HTS}}}
\newcommand\HT{{\mathrm{HT}}}
\newcommand\VS{{\mathrm{VS}}}
\newcommand\V{{\mathrm{V}}}
\newcommand\U{{\mathrm{U}}}
\newcommand\QTS{{\mathrm{QTS}}}
\newcommand\QT{{\mathrm{QT}}}
\newcommand\VHS{{\mathrm{VHS}}}
\newcommand\VHP{{\mathrm{VHP}}}
\newcommand\VH{{\mathrm{VH}}}
\newcommand\UU{{\mathrm{UU}}}
\newcommand\OS{{\mathrm{OS}}}
\newcommand\OD{{\mathrm{O}}}
\newcommand\OOS{{\mathrm{OOS}}}
\newcommand\OO{{\mathrm{OO}}}
\newcommand\UOS{{\mathrm{UOS}}}
\newcommand\UO{{\mathrm{UO}}}
\newcommand\VO{{\mathrm{VO}}}
\newcommand\VOS{{\mathrm{VOS}}}
\newcommand\DS{{\mathrm{DS}}}
\newcommand\GL{\mathbf{GL}}    % general linear group
\newcommand\Symp{\mathbf{Sp}}  % symplectic group
\newcommand\Orth{\mathbf{O}}   % orthogonal group
\newcommand\Sym{\mathcal{S}}   % symmetrci group

\title{
Enumeration of Symmetry Classes of Alternating Sign Matrices
 and Characters of Classical Groups
}

\author{
Soichi OKADA\thanks{
Graduate School of Mathematics, Nagoya University,
e-mail: okada@math.nagoya-u.ac.jp
}
}

\date{}

\begin{document}

\maketitle

\begin{abstract}
An alternating sign matrix is a square matrix with entries $1$, $0$ and $-1$
 such that the sum of the entries in each row and each column is equal to $1$
 and the nonzero entries alternate in sign along each row and each column.
To some of the symmetry classes of alternating sign matrices
 and their variations, G.~Kuperberg associate square ice models
 with appropriate boundary conditions,
 and give determinanat and Pfaffian formulae for the partition functions.
In this paper, we utilize several determinant and Pfaffian identities
 to evaluate Kuperberg's determinants and Pfaffians,
 and express the round partition functions in terms of
 irreducible characters of classical groups.
In particular, we settle a conjecture on the number of
 vertically and horizontally symmetric alternating sign matrices (VHSASMs).
\end{abstract}

\section{
Introduction
}

An {\it alternating sign matrix} (or ASM for short) is a square matrix
 satisfying the following three conditions :
\begin{roster}{(c)}
\item[(a)]
All entries are $1$, $-1$ or $0$.
\item[(b)]
Every row and column have sum $1$.
\item[(c)]
In every row and column, the nonzero entries alternate in sign.
\end{roster}
Let $\ASM_n$ be the set of $n \times n$ ASMs.
This notion of alternating sign matrices was introduced by
 D.~Robbins and H.~Rumsey \cite{RR} in a study of Dodgson's condensation
 formula for evaluating determinants.
W.~Mills, Robbins and Rumsey \cite{MRR} conjectured a formula of the number of
 $n \times n$ alternating sign matrices.
After more than 10 years, this conjecture was settled by
 D.~Zeilberger \cite{Z} and G.~Kuperberg \cite{K1} in completely different
 ways.
(See \cite{Br} for the history of ASMs and related topics.)

\begin{theorem} (Zeilberger \cite{Z}, Kuperberg \cite{K1})
The number of $n \times n$ ASMs is given by
$$
\# \ASM_n = \prod_{k=0}^{n-1} \frac{(3k+1)!}{(n+k)!}.
$$ 
\end{theorem}

Kuperberg's proof is based on a bijection between ASMs and square-ice states
 in the 6-vertex model with domain wall boundary condition,
 and on the Izergin-Korepin determinant formula for the partition function
 of this model.

The dihedral group $D_8$ of order $8$ acts on the set $\ASM_n$ of all ASMs
 as symmetries of the square.
Each subgroup of $D_8$ gives rise to symmetry classes of ASMs.
There are $7$ conjugacy classes of nontrivial subgroups of $D_8$
 and it is enough to consider the following symmetry classes of ASMs.
\begin{roster}{(WWW)}
\item[(HTS)]
{\it Half-turn symmetric} ASMs (HTSASMs)
 that are invariant under a $180^\circ$ rotation.
\item[(QTS)]
{\it Quarter-turn symmetric} ASMs (QTSASMs)
 that are invariant under a $90^\circ$ rotation.
\item[(VS)]
{\it Vertically symmetric} ASMs (VSASMs)
 that are invariant under a flip around the vertical axis.
\item[(VHS)]
{\it Vertically and horizontally symmetric} ASMs (VHSASMs)
 that are invariant under flips around both the vertical axis
 and the horizontal axis.
\item[(DS)]
{\it Diagonally symmetric} ASMs (DSASMs)
 that are symmetric in the main diagonal.
\item[(DAS)]
{\it Diagonally and antidiagonally} symmetric ASMs (DASASMs)
 that are symmetric in both diagonals.
\item[(TS)]
{\it Totally symmetric} ASMs (TSASMs)
 that are invariant under the full symmetry group $D_8$.
\end{roster}
For each symmetry class $\circledast = \HTS$, $\QTS$, $\cdots$,
 let $\ASM^\circledast_n$ denote the set of $n \times n$ ASMs with
 symmetry $\circledast$.

Kuperberg \cite{K2} extends his argument in \cite{K1}
 to several classes of ASMs (or their variations)
 including even-order HTSASMs, even-order QTSASMSs, VSASMs and VHSASMs.
He finds determinant and Pfaffian formulae for the partition functions
 of the square ice models corresponding to these classes of ASMs.
Also, by $q$-specialization, he evaluates determinants and Pfaffians
 and proves closed product formulae for the number of ASMs in
 many of these classes.
However, in the enumeration of VHSASMs, he only gives determinant
 formulae for the partition functions and do not succeed in proving
 the product formula conjectured by Mills.

In this article, we evaluate the Kuperberg's determinants and Pfaffians
 by applying determinant and Pfaffian identities
 involving Vandermonde-type determinants
 (see Theorem~3.3 and 3.4 in Section~3),
 some of which appeared in \cite{O} and were used
 for a study of rectangular-shaped representations of classical groups.
Then we can show that the partition functions corresponding to
 the round $1$-, $2$-, and $3$-enumerations are expressed in terms of
 irreducible characters of classical groups up to simple factors.
In particular, we obtain the following formulae for the number of
 some symmetry classes of ASMs.

\begin{theorem}
\begin{roster}{(5)}
\item[(1)]
The number of $n \times n$ ASMs is given by
$$
\# \ASM_n =
3^{-n(n-1)/2} \dim \GL_{2n}(\delta(n-1,n-1)).
$$
\item[(2)]
The number of $2n \times 2n$ HTSASMs is given by
$$
\# \ASM^{\HTS}_{2n} =
3^{-n(n-1)/2}\dim \GL_{2n}(\delta(n-1,n-1))
\cdot
3^{-n(n-1)/2} \dim \GL_{2n}(\delta(n,n-1)).
$$
\item[(3)]
The number of $4n \times 4n$ QTSASMs is given by
$$
\# \ASM^{\QTS}_{4n} =
\left( 3^{-n(n-1)/2} \dim \GL_{2n}(\delta(n-1,n-1)) \right)^3
\cdot
3^{-n(n-1)/2} \dim \GL_{2n}(\delta(n,n-1)).
$$
\item[(4)]
The number of $(2n+1) \times (2n+1)$ VSASMs is given by
$$
\# \ASM^{\VS}_{2n+1} =
3^{-n(n-1)} \dim \Symp_{4n}( \delta(n-1,n-1) ).
$$
\item[(5)]
The number of $(4n+1) \times (4n+1)$ VHSASMs is given by
$$
\# \ASM^{\VHS}_{4n+1} =
3^{-n(n-1)} \dim \Symp_{4n}( \delta(n-1,n-1) )
\cdot
2^{-2n} 3^{-n^2} \dim \tilde{\Orth}_{4n}( \delta(n+1/2,n-1/2) ).
$$
\item[(6)]
The number of $(4n+3) \times (4n+3)$ VHSASMs is given by
$$
\# \ASM^{\VHS}_{4n+3} =
3^{-n(n-1)} \dim \Symp_{4n}( \delta(n-1,n-1) )
\cdot
3^{-n^2} \dim \Symp_{4n+2}( \delta(n,n-1) ).
$$
\end{roster}
Here $\dim \GL_N(\lambda)$ (resp. $\dim \Symp_N(\lambda)$,
 $\dim \tilde{\Orth}_N(\lambda)$) denotes the dimension of
 the irreducible representation
 of $\GL_N$ (resp. $\Symp_N$, $\tilde{\Orth}_N$) with
 ``highest weight'' $\lambda$
 (see \S 2 for a precise definition) and
\begin{gather*}
\delta(n-1,n-1)
 = (n-1,n-1,n-2,n-2,\cdots,2,2,1,1),
\\
\delta(n,n-1)
 = (n,n-1,n-1,n-2,\cdots,3,2,2,1),
\\
\delta(n+1/2,n-1/2)
 = (n+1/2,n-1/2,n-1/2,n-3/2,\cdots, 5/2,3/2,3/2,1/2).
\end{gather*}
\end{theorem}

Each identity in this theorem, together with the Weyl's dimension formula,
 gives a closed product formula for for the number of the symmetry class
 of ASMs.
In particular, we settle the Mills' conjecture on the number of VHSASMs
 \cite[Section~4.2]{R}.

Also we obtain the following formulae for other classes of ASMs.
(See Section~2 for a definition of each class.)

\begin{theorem}
\begin{roster}{(3)}
\item[(1)]
The number of $2n \times 2n$ OSASMs is given by
$$
\# \ASM^{\OS}_{2n} =
3^{-n(n-1)} \dim \Symp_{4n}( \delta(n-1,n-1) ).
$$
\item[(2)]
The number of $(8n+1) \times (8n+1)$ VOSASMs is given by
$$
\# \ASM^{\VOS}_{8n+1} =
\left( 3^{-n(n-1)} \dim \Symp_{4n}( \delta(n-1,n-1) ) \right)^3
\cdot
2^{-2n} 3^{-n^2} \dim \tilde{\Orth}_{4n}( \delta(n+1/2,n-1/2) ).
$$
\item[(3)]
The number of $(8n+3) \times (8n+3)$ VOSASMs is given by
$$
\# \ASM^{\VOS}_{8n+3} =
\left( 3^{-n(n-1)} \dim \Symp_{4n}( \delta(n-1,n-1) ) \right)^3
\cdot
3^{-n^2} \dim \Symp_{4n+2}( \delta(n,n-1) ).
$$
\item[(4)]
The number of UASMs of order $2n$ is given by
$$
\# \ASM^{\U}_{2n} =
2^n 3^{-n(n-1)} \dim \Symp_{4n}( \delta(n-1,n-1) ).
$$
\item[(5)]
The number of UUASMs of order $4n$ is given by
$$
\# \ASM^{\UU}_{4n} =
3^{-n(n-1)} \dim \Symp_{4n}( \delta(n-1,n-1) )
 \cdot
3^{-n(n-1)} \dim \tilde{\Orth}_{4n+1}( \delta(n,n-1) ).
$$
\item[(6)]
The number of VHPASMs of order $4n+2$ is given by
$$
\# \ASM^{\VHP}_{4n+2} =
\left( 3^{-n(n-1)} \dim \Symp_{4n}( \delta(n-1,n-1) ) \right)^2.
$$
\item[(7)]
The number of UOSASMs of order $8n$ is given by
$$
\# \ASM^{\UOS}_{8n} =
\left( 3^{-n(n-1)} \dim \Symp_{4n}( \delta(n-1,n-1) ) \right)^3
 \cdot
3^{-n(n-1)} \dim \tilde{\Orth}_{4n+1}( \delta(n,n-1) ).
$$
\end{roster}
\end{theorem}

This theorem leads us to closed product formulae
 for the numbers of these classes of ASMs.
The product formulae in the VOSASM case are new, though
 the other case are studied in \cite{K2}.

This paper is organized as follows.
In Section~2, we review results in \cite{K2} on the partition functions
 of square ice models associated to various classes of ASMs,
 and state our main results which relate these partition functions
 with characters of classical groups.
As key tools in evaluating the determinants and Pfaffians appearing
 in the partition functions, we use determinant and Pfaffian identities
 involving the Vandermonde-type determinants, which are presented in
 Section~3.
In Section~4, we prove the main results.

\section{
Partition functions and classical group characters
}

In this section, we review results on the partition functions in \cite{K2}
 and give formulae which relate these partition functions to
 the classical group characters.

In addition to the symmetry classes of square ASMS, we consider
 the following classes of square ASMs.
\begin{roster}{(WWW)}
\item[(OS)]
{\it Off-diagonally symmetric} ASMs (OSASMs),
 that are diagonally symmetric ASMs with zeros on the main diagonal.
\item[(OOS)]
{\it Off-diagonally and off-antidiagonally symmetric} ASMs (OOSASMs),
 that are diagonally and antidiagonally symmetric ASMs with
 zeros on the main diagonal and the antidiagonal.
\item[(VOS)]
{\it Vertically and off-diagonally symmetric} ASMs (VOSASMs),
 that are vertically symmetric and diagonally symmetric
 with zeros on the main diagonal except for the center.
\end{roster}
The last class (VOSASMs) is not considered in \cite{K2},
 but this arises naturally from UOSASMs defined below.
It is clear that VOSASMs are TSASMs.
And one can show that there are no VOSASMs of order $8n+5$ or $8n+7$.

A vector $a = (a_1, \cdots, a_n)$ consisting $1$s, $0$s and $-1$s
 is an {\it alternating sign vector} if the sum of the entries is equal to $1$
 and the nonzero entries alternate in sign.
Kuperberg \cite{K2} introduces the following variations of ASMs.
\begin{roster}{(WWW)}
\item[(U)]
An {\it alternating sign matrix with U-turn boundary}
 (UASM) of order $2n$ is a $2n \times n$ matrix
 $A = (a_{ij})_{1 \le i \le 2n, 1 \le j \le n}$ satisfying
 the following conditions :
\begin{roster}{(2)}
\item[(1)]
Each column vector is an alternating sign vector.
\item[(2)]
For each $k$, the vector $(a_{2k-1,1}, a_{2k-1,2}, \cdots, a_{2k-1,n},
 a_{2k,n}, \cdots, a_{2k,2}, a_{2k,1})$ is an alternating sign vector.
\end{roster}
Let $\ASM^{\U}_{2n}$ be the set of UASMs of order $2n$.
\item[(UU)]
A {\it alternating sign matrix with U-U-turn boundary}
 (UUASM) of order $4n$ is a $2n \times 2n$ matrix
 $A = (a_{ij})_{1 \le i, j \le 2n}$ satisfying
 the following conditions :
\begin{roster}{(2)}
\item[(1)]
For each $k$, the vector $(a_{2k-1,1}, a_{2k-1,2}, \cdots, a_{2k-1,2n},
 a_{2k,2n}, \cdots, a_{2k,2}, a_{2k,1})$ is an alternating sign vector.
\item[(2)]
For each $k$, the vector $(a_{1, 2k-1}, a_{2,2k-1}, \cdots, a_{2n,2k-1},
 a_{2n,2k}, \cdots, a_{2,2k}, a_{1,2k})$ is an alternating sign vector.
\end{roster}
Let $\ASM^{\UU}_{4n}$ be the set of UUASMs of order $4n$.
\item[(VHP)]
A {\it vertically and horizontally perverse alternating sign matrix}
 (VHPASM) of order $4n+2$ is a UUASM
 $A = (a_{ij})_{1 \le i, j \le 2n}$ of order $4n$ such that
$\sum_{j=1}^{2n} a_{2k-1,j} = 0$ and $\sum_{i=1}^{2n} a_{i,2k-1} = 1$
 for $1 \le k \le n$.
Let $\ASM^{\VHP}_{4n+2}$ be the set of VHPASMs of order $4n+2$.
\item[(UOS)]
An {\it off-diagonally symmetric alternating sign matrix with U-turn boundary}
 (UOSASM) of order $8n$ is a UUASM
 $A = (a_{ij})_{1 \le i, j \le 4n}$ of order $4n$ such that
 $A$ is symmetric with zeros on the main diagonal.
Let $\ASM^{\UOS}_{8n}$ be the set of UOSASMs of order $8n$.
\end{roster}

For each class $\circledast = \HTS$, $\QTS$, $\cdots$, $\UOS$,
 we consider the generating function
$$
A^\circledast_n (x) = \sum_{A \in \ASM^\circledast_n} x^{n_\circledast(A)},
$$
where $n_\circledast(A)$ is the number of the orbits of $-1$s under symmetry,
 excluding any $-1$s that are forced by symmetry.
We are interested in the integers $A^\circledast_n (0)$,
 $A^\circledast_n (1)$, $A^\circledast_n(2)$, $A^\circledast_n(3)$,
 which are called $0$-, $1$-, $2$- and $3$-enumeration of the class
 $\circledast$ respectively.
In \cite{K2}, more parameters are introduced for some classes,
 but here we concentrate on these $x$-enumerations.

Now we give formulae for the partition functions of the square ice models
 associated to various classes of ASMs.
We use the following abbreviation :
$$
\sigma(t) = t - \frac{1}{t}.
$$

For two vectors of $n$ variables $\vectx = (x_1, \cdots, x_n)$,
 $\vecty = (y_1, \cdots, y_n)$ and parameters $a$, $b$, $c$,
 we introduce the following $n \times n$ matrices :
\begin{align*}
M(n ; \vectx, \vecty ; a)_{i,j}
 &=
\frac{ 1 }
     { \sigma(a x_i/y_j) \sigma(a y_j/x_i) },
\\
M_{\HT}(n ; \vectx, \vecty ; a)_{i,j}
 &=
\frac{ 1 }
     { \sigma(a x_i/y_j) }
+
\frac{1}
     { \sigma(a y_j/x_i) },
\\
M_{\U}(n ; \vectx, \vecty ; a)_{i,j}
 &=
\frac{ 1 }
     { \sigma(a x_i/y_j) \sigma(a y_j/x_i) }
-
\frac{ 1 }
     { \sigma(a x_i y_j) \sigma(1 /x_i y_j) },
\\
M_{\UU}(n ; \vectx, \vecty ; a,b,c)_{i,j}
 &=
\frac{ \sigma ( b / y_j )
       \sigma ( c x_i ) }
     { \sigma ( a x_i / y_j ) }
-
\frac{ \sigma ( b / y_j )
       \sigma ( c / x_i ) } 
     { \sigma ( a / x_i y_j ) }
-
\frac{ \sigma (b y_j)
       \sigma (c x_i) }
     { \sigma (a x_i y_j) }
+
\frac{ \sigma ( b y_j )
       \sigma ( c / x_i ) }
     { \sigma ( a y_j / x_i )}.
\end{align*}
We put
\begin{align*}
F(n ; \vectx, \vecty ; a)
&=
\frac{ \prod_{i,j=1}^n \sigma(a x_i/y_j) \sigma(a y_j/x_i) }
     { \prod_{1 \le i < j \le n} \sigma(x_j/x_i) \sigma(y_i/y_j) },
\\
F_V(n ; \vectx, \vecty ; a)
&=
\frac{ \prod_{i,j=1}^n \sigma(ax_i/y_j) \sigma(ay_j/x_i)
                       \sigma(ax_iy_j) \sigma(a/x_iy_j) }
     { \prod_{1 \le i < j \le n} \sigma(x_j/x_i) \sigma(y_i/y_j)
       \prod_{1 \le i \le j \le n} \sigma(1/x_ix_j) \sigma(y_iy_j) },
\end{align*}
and define the partition functions as follows :
\begin{align*}
A(n ; \vectx, \vecty ; a)
 &=
\sigma(a)^{-n^2+n}
F(n ; \vectx, \vecty ; a)
\det M(n ; \vectx, \vecty ; a),
\\
A_{\HT}^{(2)}(2n ; \vectx, \vecty ; a)
 &=
\sigma(a)^{-n^2}
F(n ; \vectx, \vecty ; a)
\det M_{\HT}(n ; \vectx, \vecty ; a),
\\
A_{\V}(2n+1 ; \vectx, \vecty ; a)
 &=
\sigma(a)^{-2n^2+2n}
F_{\V}(n ; \vectx, \vecty ; a)
\det M_{\U}(n ; \vectx, \vecty ; a),
\\
A_{\UU}^{(2)}(4n ; \vectx, \vecty ; a, b, c)
 &=
\sigma(a)^{-2n^2-n} \sigma(b/a)^{-n} \sigma(c/a)^{-n} \sigma(a^2)^{2n}
\\
&\quad\times
F_{\V}(n ; \vectx, \vecty ; a)
\det M_{\UU}(n ; \vectx, \vecty ; a, b, c),
\\
A_{\VH}^{(2)}(4n+1 ; \vectx, \vecty ; a)
 &=
\sigma(a)^{-2n^2-n}
F_{\V}(n ; \vectx, \vecty ; a)
\det M_{\UU}(n ; \vectx, \vecty ; a, a, a),
\\
A_{\VH}^{(2)}(4n+3 ; \vectx, \vecty ; a)
 &=
\sigma(a)^{-2n^2-n}
F_{\V}(n ; \vectx, \vecty ; a)
\det M_{\UU}(n ; \vectx, \vecty ; a, a^{-1}, a^{-1}),
\\
A_{\VHP}^{(2)}(4n+2 ; \vectx, \vecty ; a)
 &=
\sigma(a)^{-2n^2-n}
F_{\V}(n ; \vectx, \vecty ; a)
(-1)^n \det M_{\UU}(n ; \vectx, \vecty ; a, a, a^{-1}).
\end{align*}
We call them the determinant partition functions.

\begin{remark}
The partition functions of the square-ice models associated to
 UASMs and UUASMs computed in \cite{K2} have extra factors
$$
\prod_{i=1}^n
 \frac{\sigma(a^2 x_i^2)}{\sigma(a^2)}
 \frac{\sigma(b/y_i)}{\sigma(b)},
\quad\text{and}\quad
\prod_{i=1}^n
 \frac{\sigma(a^2 x_i^2)}{\sigma(a^2)}
 \frac{\sigma(a^2 / y_i^2)}{\sigma(a^2)},
$$
which do not affect the $x$-enumeration.
So we omit these factors in the definition of $A_{\V}$ and $A_{\UU}^{(2)}$.
\end{remark}

For a vector of $2n$ variables $\vectx = (x_1, \cdots, x_{2n})$
 and parameters $a$, $b$, $c$, we introduce the $2n \times 2n$ skew-symmetric
 matrices :
\begin{align*}
M_{\QT}^{(k)}(n ; \vectx ; a)_{ij}
 &=
\frac{ \sigma( x_j^k / x_i^k ) }
     { \sigma( a x_j / x_i ) \sigma( a x_i / x_j ) },
\\
M_{\OD}(n ; \vectx ; a)_{ij}
 &=
\frac{ \sigma(x_j/x_i) }
     { \sigma(a x_i x_j) \sigma (a / x_i x_j) },
\\
M_{\OO}(n ; \vectx ; a, b, c)_{ij}
 &=
\sigma(x_j/x_i)
\left(
 \frac{ c^2 }
      { \sigma(a x_i x_j) }
+
 \frac{ b^2 }
      { \sigma(a / x_i x_j) }
\right),
\\
M_{\UO}^{(1)}(n ; \vectx ; a)_{i,j}
&=
\sigma(x_j / x_i) \sigma(x_i x_j)
\left(
 \frac{1}{ \sigma(a x_i x_j) \sigma(a / x_i x_j) }
-
 \frac{1}{ \sigma(a x_j / x_i) \sigma(a x_i / x_j) }
\right),
\\
M_{\UO}^{(2)}(n ; \vectx ; a, c)_{i,j}
&=
\sigma(x_j / x_i) \sigma(x_i x_j)
\\
&\quad\times \left(
 \frac{ \sigma(c x_i) \sigma(c x_j) }
      { \sigma(a x_i x_j) }
-
 \frac{ \sigma(c x_i) \sigma(c / x_j) }
      { \sigma(a x_i / x_j) }
-
 \frac{ \sigma(c / x_i) \sigma(c x_j) }
      { \sigma(a x_j / x_j) }
+
 \frac{ \sigma(c / x_i) \sigma(c / x_j) }
      { \sigma(a / x_i x_j) }
\right).
\end{align*}
We put
\begin{align*}
F_{\QT}(n ; \vectx ; a)
&=
\frac{ \prod_{1 \le i < j \le 2n} \sigma(a x_j / x_i) \sigma(a x_i / x_j) }
     { \prod_{1 \le i < j \le 2n} \sigma(x_j / x_i) },
\\
F_{\OD}(n ; \vectx ; a)
&=
\frac{ \prod_{1 \le i < j \le 2n}
         \sigma(a x_i x_j) \sigma (a / x_i x_j)
     }
     { \prod_{1 \le i < j \le 2n}
         \sigma(x_j/x_i)
     }
\\
F_{\UO}(n ; \vectx ; a)
&=
\frac{ \prod_{1 \le i < j \le 2n}
         \sigma(ax_i/x_i) \sigma(ax_j/x_i)
         \sigma(ax_ix_j) \sigma(a/x_ix_j)
     }
     { \prod_{1 \le i < j \le 2n} \sigma(x_j/x_i)
       \prod_{1 \le i \le j \le 2n} \sigma(x_i x_j)
     },
\end{align*}
and define the Pfaffian partition functions as follows :
\begin{align*}
A_{\QT}^{(k)}(4n ; \vectx ; a)
 &=
\sigma(a)^{-2n^2+2n}
F_{\QT}(n ; \vectx ; a)
\Pf M_{\QT}^{(k)}(n ; \vectx ; a),
\\
A_{\OD}(2n ; \vectx ; a)
 &=
\sigma(a)^{-2n^2+2n}
F_{\OD}(n ; \vectx ; a)
\Pf M_{\OD}(n ; \vectx ; a),
\\
A_{\OO}^{(2)}(4n ; \vectx ; a, b, c)
&=
c^{-2n} \sigma(a)^{-2n^2+n}
F_{\OD}(n ; \vectx ; a)
\Pf M_{\OO}(n ; \vectx ; a, b, c),
\\
A_{\UO}^{(1)}(8n ; \vectx ; a)
&=
\sigma(a)^{-4n^2+4n}
F_{\UO}(n ; \vectx ; a)
\Pf M_{\UO}^{(1)}(n ; \vectx ; a),
\\
A_{\UO}^{(2)}(8n ; \vectx ; a, c)
&=
\sigma(a)^{-4n^2+n} \sigma(c/a)^{-2n} \sigma(a^2)^{2n}
F_{\UO}(n ; \vectx ; a)
\Pf M_{\UO}^{(2)}(n ; \vectx ; a, c),
\\
A_{\VO}^{(2)}(8n+1 ; \vectx ; a)
&=
\sigma(a)^{-4n^2+n}
F_{\UO}(n ; \vectx ; a)
\Pf M_{\UO}^{(2)}(n ; \vectx ; a, a),
\\
A_{\VO}^{(2)}(8n+3 ; \vectx ; a)
&=
\sigma(a)^{-4n^2+n}
F_{\UO}(n ; \vectx ; a)
\Pf M_{\UO}^{(2)}(n ; \vectx ; a, a^{-1}).
\end{align*}
We call them the Pfaffian partition functions.

Kuperberg \cite{K2} proves that the $x$-enumerations are obtained
 from these partition functions by specializing all the spectral parameters
 $x_1, \cdots, x_n, y_1, \cdots, y_n$ (or $x_1, \cdots, x_{2n}$)
 to $1$.

\begin{theorem} (Kuperberg)
Let $\vectone = (1, 1, \cdots, 1)$.
If $x = a^2 + 2 + a^{-2}$, then we have
\begin{align*}
A_n (x)
 &= A(n ; \vectone, \vectone ; a),
\\
A^{\HTS}_{2n} (x)
 &=
A(n ; \vectone, \vectone ; a)
A_{\HT}^{(2)}(2n ; \vectone, \vectone ; a),
\\
A_{2n+1}^{\VS} (x)
 &=
 A_{\V}(2n+1 ; \vectone, \vectone ; a),
\\
A_{4n+1}^{\VHS} (x)
 &=
 A_{\V}(2n+1 ; \vectone, \vectone ; a)
 A_{\VH}^{(2)}(4n+1 ; \vectone, \vectone ; a),
\\
A_{4n+3}^{\VHS} (x)
 &=
 A_{\V}(2n+1 ; \vectone, \vectone ; a)
 A_{\VH}^{(2)}(4n+3 ; \vectone, \vectone ; a),
\\
A_{2n}^{\U} (x)
 &=
 2^n 
 A_{\V}(2n+1 ; \vectone, \vectone ; a),
\\
A^{\UU}_{4n}(x)
 &=
 A_{\V}(2n+1 ; \vectone, \vectone : a)
 A_{\UU}^{(2)}(4n ; \vectone, \vectone ; a, \sqrt{-1}, \sqrt{-1}),
\\
A_{4n+2}^{\VHP} (x)
 &=
 A_{\V}(2n+1 ; \vectone, \vectone ; a)
 A_{\VHP}^{(2)}(4n+2 ; \vectone, \vectone ; a),
\end{align*}
Also we have
\begin{align*}
A^{\QTS}_{4n}(x)
 &=
 A_{\QT}^{(1)}(4n ; \vectone ; a)
 A_{\QT}^{(2)}(4n ; \vectone ; a),
\\
A^{\OS}_{2n}(x)
 &=
 A_{\OD}(2n ; \vectone ; a),
\\
A^{\OOS}_{4n}(x)
 &=
 A^{\OD}(2n ; \vectone ; a)
 A_{\OO}^{(2)}(4n ; \vectone ; a, b, b),
\\
A^{\UOS}_{8n}(x)
 &=
 A_{\UO}^{(1)}(8n ; \vectone ; a)
 A_{\UO}^{(2)}(8n ; \vectone ; a, \sqrt{-1}),
\\
A^{\VOS}_{8n+1}(x)
 &=
 A_{\UO}^{(1)}(8n ; \vectone ; a)
 A_{\VO}^{(2)}(8n+1 ; \vectone ; a),
\\
A^{\VOS}_{8n+3}(x)
 &=
 A_{\UO}^{(1)}(8n ; \vectone ; a)
 A_{\VO}^{(2)}(8n+3 ; \vectone ; a).
\end{align*}
\end{theorem}

\begin{remark}
The last two identities for VOSASMs were not treated in \cite{K2},
 but can be proved in a way similar to the proof for VHSASMs.
\end{remark}

Let $\zeta_n = \exp (2 \pi \sqrt{-1} / n)$ be the primitive $n$th
 root of unity.
Then the correspondence between $a$ and $x = a^2 + 2 + a^{-2}$ are
 given as follows:
$$
\begin{array}{c|c|c|c|c}
a & \zeta_4 & \zeta_6 & \zeta_8 & \zeta_{12} \\
\hline
x = a^2 + 2 + a^{-2} & 0 & 1 & 2 & 3
\end{array}
$$

To state our results, we introduce the irreducible characters of
 classical groups.
A partition (resp. half-partition) is a non-increasing sequence
 $\lambda = (\lambda_1, \cdots, \lambda_n)$ of non-negative integers
 $\lambda_i \in \Int$ (resp. non-negative half-integers $\lambda_i
 \in \Int+1/2$).
If $\lambda$ and $\mu$ are both partitions (resp. both half-partitions),
 then $\lambda \cup \mu$ denote the partition (resp. half-partition)
 obtained by rearranging the entries of $\lambda$ and $\mu$
 in decreasing order.
For a integer (or a half-integer) $r$ and a positive integer $k$,
 we define a partition (or a half-partion) $\delta(r)$ and $\delta^2(r)$ by
 putting
$$
\delta(r) = (r, r-1, r-2, r-3, \cdots),
\quad
\delta^2(r) = (r, r-2, r-4, r-6, \cdots).
$$
Also we define
$$
\delta(r,s) = \delta(r) \cup \delta(s),
\quad
\delta^2(r,s) = \delta^2(r) \cup \delta^2(s).
$$

For a sequence $\alpha = (\alpha_1, \cdots, \alpha_n)$ of
 integers or half integers,
 and a vector $\vectx = (x_1, \cdots, x_n)$ of $n$ indeterminates,
we define $n \times n$ matrices $V(\alpha ; \vectx)$ and
 $W_{\pm}(\alpha ; \vectx)$ by putting
$$
V(\alpha ;\vectx)
 = \left( x_i^{\alpha_j} \right)_{1 \le i, j \le n},
\quad
W^{\pm}(\alpha ; \vectx)
 = \left( x_i^{\alpha_j} \pm x_i^{-\alpha_j} \right)_{1 \le i, j \le n}.
$$
For a partition $\lambda$ with length $\le n$, we define
\begin{align}
\GL_n(\lambda;\vectx)
 &=
 \frac{ \det V(\lambda + \delta(n-1) ; \vectx) }
      { \det V(\delta(n-1) ; \vectx) },
\label{GL}
\\
\Symp_{2n}(\lambda;\vectx)
 &=
 \frac{ \det W^-(\lambda + \delta(n) ; \vectx) }
      { \det W^-(\delta(n) ; \vectx) }.
\label{Sp}
\end{align}
Then $\GL_n(\lambda;\vectx)$ (resp. $\Symp_{2n}(\lambda;\vectx)$)
 is the character of the irreducible representation $\GL_n(\lambda)$
 (resp. $\Symp_{2n}(\lambda)$) of the general linear group $\GL_n(\Comp)$
 (resp. the symplectic group $\Symp_{2n}(\Comp)$).
If $\lambda$ is a partition with length $\le n$ or
 a half-partition of length $n$, we define
\begin{align}
\tilde{\Orth}_{2n+1}(\lambda;\vectx)
 &=
 \frac{ \det W^-(\lambda + \delta(n-1/2) ; \vectx) }
      { \det W^-(\delta(n-1/2) ; \vectx) },
\label{O1}
\\
\tilde{\Orth}_{2n}(\lambda;\vectx)
 &=
 \begin{cases}
 \dfrac{ \det W^+(\lambda + \delta(n-1) ; \vectx) }
       { (1/2) \det W^+(\delta(n-1) ; \vectx) }
 &\text{if $\lambda_n \neq 0$,} \\
 \dfrac{ \det W^+(\lambda + \delta(n-1) ; \vectx) }
       { \det W^+(\delta(n-1) ; \vectx) }
 &\text{if $\lambda_n = 0$.} \\
 \end{cases}
\label{O2}
\end{align}
Then $\tilde{\Orth}_N(\lambda;\vectx)$ is the character of the irreducible
 representation $\tilde{\Orth}_N(\lambda)$ of the double cover
 $\tilde{\Orth}_N$ of the orthogonal group $\Orth_N$.

Now we are in position to state our main results.
Theorem~1.2 and 1.3 in Introduction immediately follows
 from the following Theorem and Theorem~2.1.
For a vector $\vectx = (x_1, \cdots, x_n)$ and a integer $k$,
 we put
$$
\vectx^k = (x_1^k, \cdots, x_n^k).
$$
First we give formulae for the determinant partition functions.

\begin{theorem}
\begin{roster}{(7)}
\item[(1)]
For the partition function associated to ASMs, we have
\begin{align*}
A(n ; \vectx, \vecty ; \zeta_4)
&=
2^{-n^2+n}
 \prod_{i=1}^n x_i^{-n+1} y_i^{-n+1}
 \prod_{i,j=1}^n (x_i^2 + y_j^2)
 \cdot \perm \left( \frac{1}{x_i^2 + y_j^2} \right)_{1 \le i, j \le n},
\\
A(n ; \vectx, \vecty ; \zeta_6)
&=
3^{-n(n-1)/2}
 \prod_{i=1}^n x_i^{-n+1} y_i^{-n+1}
 \GL_{2n}(\delta(n-1,n-1) : \vectx^2, \vecty^2),
\\
A(n ; \vectx, \vecty ; \zeta_8)
&=
2^{-n(n-1)/2}
 \prod_{i=1}^n x_i^{-n+1} y_i^{-n+1}
 \prod_{1 \le i < j \le n} (x_i^2 + x_j^2) (y_i^2 + y_j^2),
\\
A(n ; \vectx, \vectx ; \zeta_{12})
&=
\prod_{i=1}^n x_i^{-2n+2}
 \GL_n(\delta(p,p-1) ; \vectx^4) \GL_n(\delta(q,q) ; \vectx^4),
\end{align*}
where $\perm$ denotes the permanent of a matrix (see (\ref{perm})),
and $p$ and $q$ are the largest integers not exceeding $n/2$
 and $(n-1)/2$ respectively.
\item[(2)]
For the partition function associated to HTSASMs of order $2n$, we have
\begin{align*}
A_{\HT}^{(2)}(2n ; \vectx, \vecty ; \zeta_4)
&=
2^{-n^2+n}
\prod_{i=1}^n x_i^{-n} y_i^{-n}
\prod_{i,j=1}^n (x_i^2 + y_j^2),
\\
A_{\HT}^{(2)}(2n ; \vectx, \vecty ; \zeta_6)
&=
3^{-n(n-1)/2}
 \prod_{i=1}^n x_i^{-n} y_i^{-n}
 \GL_{2n}(\delta(n, n-1) ; \vectx^2, \vecty^2),
\\
A_{\HT}^{(2)}(n ; \vectx, \vectx ; \zeta_8)
&=
2^{-n(n-1)/2+n}
 \GL_n(\delta^2(n,n-2) ; \vectx) \GL_n( \delta^2(n-1,n-1) ; \vectone).
\end{align*}
\item[(3)]
For the partition function associated to VHSASMs, we have
\begin{align*}
A_{\V}(2n+1 ; \vectx, \vecty ; \zeta_4)
&=
2^{-2n^2+2n}
 \prod_{i=1}^n x_i^{-2n+2} y_i^{-2n+2}
\\
&\quad\times
 \prod_{i,j=1}^n (x_i^2 + y_j^2)(1 + x_i^2 y_j^2)
 \perm \left(
  \frac{ 1 }
       { (x_i^2 + y_j^2) (1 + x_i^2 y_j^2) }
 \right)_{1 \le i, j \le n},
\\
A_{\V}(2n+1 ; \vectx, \vecty ; \zeta_6)
&=
3^{-n(n-1)}
\Symp_{4n}( \delta(n-1,n-1) ; \vectx^2, \vecty^2 ),
\\
A_{\V}(2n+1 ; \vectx, \vecty ; \zeta_8)
&=
2^{-n(n-1)}
\prod_{i=1}^n x_i^{-2n+2} y_i^{-2n+2}
\\
&\quad\times
\prod_{1 \le i < j \le n}
 ( x_i^2 + x_j^2 ) ( 1 + x_i^2 x_j^2 )
 ( y_i^2 + y_j^2 ) ( 1 + y_i^2 y_j^2 ),
\\
A_{\V}(2n+1 ; \vectx, \vectx ; \zeta_{12})
&=
\frac{ 1 }
     { \prod_{i=1}^n (x_i^2 + x_i^{-2}) }
\\
&\quad\times
 \tilde{\Orth}_{2n+1}( \delta(n/2,n/2-1) ; \vectx^4)
 \tilde{\Orth}_{2n+1}( \delta((n-1)/2,(n-1)/2) ; \vectx^4).
\end{align*}
\item[(4)]
For the partition function associated to VHSASMs of order $4n+1$, we have
\begin{align*}
A_{\VH}^{(2)}(4n+1 ; \vectx, \vecty ; \zeta_4)
&=
2^{-2n^2}
\prod_{i=1}^n x_i^{-2n} y_i^{-2n}
\prod_{i, j=1}^n ( x_i^2 + y_j^2 ) ( 1 + x_i^2 y_j^2 ),
\\
A_{\VH}^{(2)}(4n+1 ; \vectx, \vecty ; \zeta_6)
&=
3^{-n^2}
\frac{ 1 }
     { \prod_{i=1}^n (x_i + x_i^{-1}) (y_i + y_i^{-1}) }
\tilde{\Orth}_{4n}( \delta(n+1/2,n-1/2) ; \vectx^2,\vecty^2),
\\
A_{\VH}^{(2)}(4n+1 ; \vectx, \vectx ; \zeta_8)
&=
2^{-n(n-1)}
\frac{ 1 }
     { \prod_{i=1}^n (x_i + x_i^{-1})^2 }
\\
&\quad\times
\tilde{\Orth}_{2n}( \delta^2 ( n+1/2, n-3/2 ) ; \vectx^2)
\tilde{\Orth}_{2n}( \delta^2 ( n-1/2, n-1/2 ) ; \vectx^2).
\end{align*}
\item[(5)]
For the partition function associated to VHSASMs of order $4n+3$, we have
\begin{align*}
A_{\VH}^{(2)}(4n+3 ; \vectx, \vecty ; \zeta_4)
&=
2^{-2n^2}
\prod_{i=1}^n x_i^{-2n} y_i^{-2n}
\prod_{i, j=1}^n ( x_i^2 + y_j^2 ) ( 1 + x_i^2 y_j^2 ),
\\
A_{\VH}^{(2)}(4n+3 ; \vectx, \vecty ; \zeta_6)
&=
3^{-n^2}
\Symp_{4n+2}( \delta(n,n-1) ; \vectx^2, \vecty^2, 1),
\\
A_{\VH}^{(2)}(4n+3 ; \vectx, \vectx ; \zeta_8)
&=
2^{-n^2+n}
\tilde{\Orth}_{2n}( \delta^2 ( n, n-2 ) ; \vectx^2)
\tilde{\Orth}_{2n}( \delta^2 ( n-1, n-1 ) ; \vectx^2).
\end{align*}
\item[(6)]
For the partition function associated to UUASMs of order $4n$, we have
\begin{align*}
A_{\UU}^{(2)}(4n ; \vectx, \vecty ; \zeta_4, \zeta_4, \zeta_4)
&=
2^{-2n^2+2n} \prod_{i=1}^n x_i^{-2n} y_i^{-2n}
\prod_{i,j=1}^n (x_i^2 + y_j^2)(1 + x_i^2 y_j^2),
\\
A_{\UU}^{(2)}(4n ; \vectx, \vecty ; \zeta_6, \zeta_4, \zeta_4)
&=
3^{-n^2+n}
\tilde{\Orth}_{4n+1}( \delta(n,n-1) ; \vectx^2, \vecty^2 ),
\\
A_{\UU}^{(2)}(4n ; \vectx, \vectx ; \zeta_8, \zeta_4, \zeta_4)
&=
2^{-n^2+2n}
\tilde{\Orth}_{2n+1}( \delta^2(n,n-2) ; \vectx^2)
\tilde{\Orth}_{2n+1}( \delta^2(n-1,n-1) ; \vectx^2).
\end{align*}
\item[(7)]
For the partition function associated to VHPASMs, we have
\begin{align*}
A_{\VHP}^{(2)}(4n+2 ; \vectx, \vecty ; \zeta_4)
&=
2^{-2n^2}
\prod_{i=1}^n x_i^{-2n} y_i^{-2n}
\prod_{i,j=1}^n (x_i^2 + y_j^2) (1 + x_i^2 y_j^2),
\\
A_{\VHP}^{(2)}(4n+2 ; \vectx, \vecty ; \zeta_6)
&=
3^{-n^2} \prod_{i=1}^n \left( y_i^2 + 1 + y_i^{-2} \right)
\Symp_{4n}( \delta(n-1,n-1) ; \vectx^2, \vecty^2 ),
\\
A_{\VHP}^{(2)}(4n+2 ; \vectx, \vectx ; \zeta_8)
&=
2^{-n^2+n}
\GL_{2n}( \delta^2(2n-2,2n-2) ; \vectx^2,\vectx^{-2}).
\end{align*}
\end{roster}
\end{theorem}

Next theorem gives formulae for the Pfaffian partition functions.

\begin{theorem}
\begin{roster}{(7)}
\item[(1)]
For the partition functions associated to QTASMs, we have
\begin{align*}
A_{\QT}^{(1)}(4n ; \vectx ; \zeta_4)
&=
2^{-2n^2+2n}
\prod_{i=1}^{2n} x_i^{-2n+2}
\prod_{1 \le i < j \le 2n} (x_i^2 + x_j^2)
\Hf \left( \frac{1}{x_i^2 + x_j^2} \right)_{1 \le i, j \le 2n},
\\
A_{\QT}^{(1)}(n ; \vectx ; \zeta_6)
&=
3^{-n^2+n} \prod_{i=1}^{2n} x_i^{-2n+2}
\GL_{2n}( \delta(n-1,n-1) ; \vectx^2 )^2,
\\
A_{\QT}^{(1)}(4n ; \vectx ; \zeta_8)
&=
2^{-n^2+n} \prod_{i=1}^{2n} x_i^{-2n+2}
\GL_{2n}( \delta^2(2n-2,2n-2) ; \vectx^2),
\\
A_{\QT}^{(2)}(4n ;\vectx ; \zeta_4 )
&=
2^{-2n^2+2n}
\prod_{i=1}^{2n} x_i^{-2n+1}
\prod_{1 \le i < j \le 2n} (x_i^2 + x_j^2),
\\
A_{\QT}^{(2)}(n ; \vectx ; \zeta_6)
&=
3^{-n^2+n} \prod_{i=1}^{2n} x_i^{-2n+1}
\GL_{2n}( \delta(n-1,n-1) ; \vectx^2 )
\GL_{2n}( \delta(n,n-1) ; \vectx^2),
\\
A_{\QT}^{(2)}(4n ; \vectx ; \zeta_8)
&=
2^{-n^2+n}
\prod_{i=1}^{2n} x_i^{-2n+1}
\prod_{1 \le i < j \le 2n} (x_i^2 + x_j^2),
\end{align*}
where $\Hf$ denotes the Hafnian of a symmetric matrix (see (\ref{Hf})).
\item[(2)]
For the partition function associated to OSASMs, we have
\begin{align*}
A_{\OD}(n ; \vectx ; \zeta_4)
&=
2^{-2n^2+2n}
\prod_{i=1}^n x_i^{-2n+2}
\prod_{1 \le i < j \le 2n} (1 + x_i^2 x_j^2)
\Hf \left(
 \frac{ 1 }
      { (1 + x_i^2 x_j^2)^2 }
\right)_{1 \le i, j \le 2n},
\\
A_{\OD}(2n ; \vectx ; \zeta_6)
&=
3^{-n^2+n}
\Symp_{4n}( \delta(n-1,n-1) ; \vectx^2).
\end{align*}
\item[(3)]
For the partition function associated to OOSASMs, we have
$$
A_{\OO}^{(2)}(4n ; \vectx ; \zeta_4, b, b)
=
2^{-2n^2+2n} \prod_{i=1}^{2n} x_i^{-2n+1}
\prod_{1 \le i < j \le 2n} (1 + x_i^2 x_j^2).
$$
\item[(4)]
For the partition functions associated to UOSASMs, we have
\begin{align*}
A_{\UO}^{(1)}(8n ; \vectx ; \zeta_4)
&=
2^{-4n^2+4n}
\prod_{i=1}^n x_i^{-4n+4}
\\
&\quad\times
\prod_{1 \le i < j \le 2n}
 (x_j^2 + x_i^2)(1 + x_i^2 x_j^2)
\Hf \left(
 \frac{ 1 }
      { (x_j^2 + x_i^2)(1 + x_i^2 x_j^2) }
\right)_{1 \le i, j \le 2n},
\\
A_{\UO}^{(1)}(8n ; \vectx ; \zeta_6)
&=
3^{-2n^2+2n}
\Symp_{4n}( \delta(n-1,n-1) ; \vectx^2 )^2,
\\
A_{\UO}^{(2)}(8n ; \vectx ; \zeta_4, \zeta_4)
&=
2^{-4n^2+4n} \prod_{i=1}^n x_i^{-4n+2}
\prod_{1 \le i < j \le 2n} (x_j^2 + x_i^2)(1 + x_i^2 x_j^2),
\\
A_{\UO}^{(2)}(8n ; \vectx ; \zeta_6, \zeta_4)
&=
3^{-2n^2+2n}
\Symp_{4n}( \delta(n-1,n-1) ; \vectx^2 )
\tilde{\Orth}_{4n+1}( \delta(n,n-1) ; \vectx^2 ).
\end{align*}
\item[(5)]
For the partition function associated to VOSASMs of order $8n+1$, we have
\begin{align*}
A_{\VO}^{(2)}(8n+1 ; \vectx ; \zeta_4)
&=
2^{-4n^2+2n}
\prod_{i=1}^{2n} x_i^{-4n+2}
\prod_{1 \le i < j \le 2n}
 (x_j^2 + x_i^2) (1 + x_i^2 x_j^2),
\\
A_{\VO}^{(2)}(8n+1 ; \vectx ; \zeta_6)
&=
3^{-2n(2n-1)/2}
\frac{ 1 }
     { \prod_{i=1}^{2n} (x_i + x_i^{-1}) }
\\
&\quad\times
\Symp_{4n}( \delta(n-1,n-1) ; \vectx^2)
\tilde{\Orth}_{4n}( \delta(n+1/2,n-1/2) ; \vectx^2).
\end{align*}
\item[(6)]
For the partition function associated to VOSASMs of order $8n+1$, we have
\begin{align*}
A_{\VO}^{(2)}(8n+3 ; \vectx ; \zeta_4)
&=
2^{-4n^2+2n}
\prod_{i=1}^{2n} x_i^{-4n+2}
\prod_{1 \le i < j \le 2n}
 (x_j^2 + x_i^2) (1 + x_i^2 x_j^2),
\\
A_{\VO}^{(2)}(8n+3 ; \vectx ; \zeta_6)
&=
3^{-2n^2+n}
\Symp_{4n}( \delta(n-1,n-1) ; \vectx^2 )
\Symp_{4n+2}( \delta(n,n-1) ; \vectx^2, 1).
\end{align*}
\end{roster}
\end{theorem}

By combining Theorem~2.3 and 2.4 with Theorem~2.1, we obtain formulae
 for $0$-, $1$-, $2$-, and $3$-enumerations in terms of the dimensions
 of irreducible representations of classical groups.
In particular, we obtain Theorem~1.2 and 1.3 in Introduction.

\begin{remark}
By using different techniques,
 Y.~Stroganov and A.~Razumov \cite{Str}, \cite{RS} obtained
 formulae of $A(n ; \vectx, \vecty ; \zeta_6)$,
 $A_{\V}(2n+1 ; \vectx, \vecty ; \zeta_6)$ and
 $A_{\OD}(2n ; \vectx ; \zeta_6)$ in terms of Vandermonde-type determinants,
 which immediately imply the corredponding formulae in Theorems~2.4 and 2.5.
\end{remark}

\section{
Determinant and Pfaffian identities
}

In this section, we collect determinant and Pfaffian identities,
 which will be used in the evaluation of the determinants and Pfaffians
 in the partition functions introduced in Section~2.

The determinant and permanent of a square matrix
 $A = (a_{ij})_{1 \le i, j \le n}$ are defined by
\begin{align}
\det A &=
 \sum_{\sigma \in \Sym_n}
  \sgn(\sigma) a_{1 \sigma(1)} a_{2 \sigma(2)} \cdots a_{n \sigma(n)},
\label{det}
\\
\perm A &=
 \sum_{\sigma \in \Sym_n}
  a_{1 \sigma(1)} a_{2 \sigma(2)} \cdots a_{n \sigma(n)},
\label{perm}
\end{align}
where $\Sym_n$ is the symmetric groups of degree $n$.
And the Pfaffian of a skew-symmetric matrix $A = (a_{ij})_{1 \le i, j \le 2n}$
 and the Hafnian of a symmetric matrix $B = (b_{ij})_{1 \le i, j \le 2n}$
 are given by
\begin{align}
\Pf A &=
 \sum_{\sigma \in {\mathcal{F}}_{2n}}
  \sgn(\sigma) a_{\sigma(1)\sigma(2)} a_{\sigma(3)\sigma(4)} \cdots
  a_{\sigma(2n-1)\sigma(2n)},
\label{Pf}
\\
\Hf B &=
 \sum_{\sigma \in {\mathcal{F}}_{2n}}
  b_{\sigma(1)\sigma(2)} b_{\sigma(3)\sigma(4)} \cdots
  b_{\sigma(2n-1)\sigma(2n)},
\label{Hf}
\end{align}
where $\mathcal{F}_{2n}$ is the set of all permutations $\sigma$ satisfying
 $\sigma(1) < \sigma(3) < \cdots < \sigma(2n-1)$ and
 $\sigma(2i-1) < \sigma(2i)$ for $1 \le i \le n$.

First we recall the Cauchy's determinant identity \cite{C},
 the Schur's Pfaffian identity \cite{S} and its variant
 (\cite{LLT}, \cite{St}).

\begin{lemma}
\begin{align}
\det \left(
 \frac{ 1 }
      { x_i + y_j }
 \right)_{1 \le i, j \le n}
 &=
\frac{ \prod_{1 \le i < j \le n} (x_j - x_i) (y_j - y_i) }
     { \prod_{i,j=1}^n (x_i + y_j) },
\label{C1}
\\
\Pf \left( \frac{x_j - x_i}{x_j + x_i} \right)_{1 \le i, j \le 2n}
 &=
\prod_{1 \le i < j \le 2n} \frac{x_j - x_i}{x_j+x_i},
\label{S1}
\\
\Pf \left(
 \frac{ x_j - x_i}
      { 1 - x_i x_j }
\right)_{1 \le i, j \le 2n}
 &=
\prod_{1 \le i < j \le 2n}
 \frac{ x_j - x_i }
      { 1 - x_i x_j }.
\label{S2}
\end{align}
\end{lemma}

The identities in the next lemma will be used to evaluate some of
 the determinants and Pfaffians appearing in the $0$-enumerations.
The first identity (\ref{B1}) goes back to C.~Borchardt \cite{B},
 and its Pfaffian-Hafnian analogues (\ref{I1}) and (\ref{I2}) are
 given in \cite{IKO}.

\begin{lemma}
\begin{align}
\det \left( \frac{1}{(x_i + y_j)^2} \right)_{1 \le i, j \le n}
 &=
\frac{\prod_{1 \le i < j \le n} (x_j - x_i) (y_j - y_i) }
     {\prod_{i,j=1}^n (x_i + y_j) }
 \cdot
 \perm \left( \frac{1}{x_i + y_j} \right)_{1 \le i, j \le n},
\label{B1}
\\
\Pf \left( \frac{x_j - x_i}{(x_j + x_i)^2} \right)_{1 \le i, j \le 2n}
 &=
\prod_{1 \le i < j \le 2n} \frac{x_j - x_i}{x_j+x_i}
\cdot
\Hf \left( \frac{1}{x_j + x_i} \right)_{1 \le i, j \le 2n},
\label{I1}
\\
\Pf \left(
 \frac{ x_j - x_i }
      { (1 - x_i x_j)^2 }
\right)_{1 \le i, j \le 2n}
 &=
\prod_{1 \le i < j \le 2n}
 \frac{ x_j - x_i }
      { 1 - x_i x_j }
\cdot
\Hf \left(
 \frac{ 1 }
      { 1 - x_i x_j }
\right)_{1 \le i, j \le 2n}.
\label{I2}
\end{align}
\end{lemma}

The following two theorems are the key to evaluate the determinants and
 Pfaffians appearing in the round $1$-, $2$-, and $3$-enumerations.
The identities (\ref{D1}), (\ref{D2}), (\ref{P1}) and (\ref{P2})
 already appeared in \cite{O} and their specializations are
 \cite[Theorem~16, 17]{K2}.

For $\vectx = (x_1, \cdots, x_n)$ and $\vecta = (a_1, \cdots, a_n)$,
let $V^{p,q}(\vectx ; \vecta)$ ($p+q=n$) and $W^n(\vectx ; \vecta)$ be
 the $n \times n$ matrices with $i$th row
$$
(1,x_i,x_i^2,\cdots, x_i^{p-1}, a_i, a_i x_i, \cdots, a_i x_i^{q-1}),
$$
$$
(1+a_i x_i^{n-1}, x_i + a_i x_i^{n-2}, \cdots, x_i^{n-1} + a_i)
$$
respectively.

\begin{theorem}
For $\vectx = (x_1, \cdots, x_n)$, $\vecty = (y_1, \cdots, y_n)$,
 $\vecta = (a_1, \cdots, a_n)$ and $\vectb = (b_1, \cdots, b_n)$,
 we have
\begin{align}
\det \left(
 \frac{b_j - a_i}{y_j - x_i}
\right)_{1 \le i, j \le n}
 &=
\frac{ (-1)^{n(n-1)/2} }{ \prod_{i,j=1}^n (y_j - x_i) }
 \det V^{n,n}(\vectx,\vecty ; \vecta,\vectb),
\label{D1}
\\
\det \left(
 \frac{ \det W^2(x_i,y_j;a_i,b_j) }
      { (1 - x_i y_j)(y_j - x_i) }
\right)_{1 \le i, j \le n}
 &=
\frac{ 1 }
     { \prod_{i,j=1}^n (1 - x_i y_j)(y_j - x_i) }
\det W^{2n}(\vectx, \vecty ; \vecta, \vectb),
\label{D2}
\\
\det \left(
 \frac{ \det W^3(x_i,y_j,z ; a_i, b_j, c) }
      { ( 1 - x_i y_j ) ( y_j - x_i ) }
\right)_{1 \le i, j \le n}
 &=
\frac{ (1+c)^{n-1} }
     { \prod_{i,j=1}^n ( 1 - x_i y_j ) ( y_j - x_i ) }
\det W^{2n+1}(\vectx,\vecty,z ; \vecta, \vectb,c).
\label{D3}
\end{align}
\end{theorem}

\begin{theorem}
For $\vectx = (x_1, \cdots, x_{2n})$, $\vecta = (a_1, \cdots, a_{2n})$,
 and $\vectb = (b_1, \cdots, b_{2n})$, we have
\begin{align}
&
\Pf \left(
 \frac{ (a_j - a_i) (b_j - b_i) }
      { x_j - x_i }
\right)_{1 \le i, j \le 2n}
 =
\frac{ 1 }
     { \prod_{1 \le i < j \le 2n} ( x_j - x_i ) }
\det V^{n,n}(\vectx ; \vecta)
\det V^{n,n}(\vectx ; \vectb),
\label{P1}
\\
&
\Pf \left(
 \frac{ \det W^2(x_i, x_j ; a_i, a_j) \det W^2(x_i, x_j ; b_i, b_j) }
      { (1 - x_i x_j)(x_j - x_i) }
\right)_{1 \le i, j \le 2n}
\notag
\\
&\qquad=
\frac{ 1 }
     { \prod_{1 \le i < j \le 2n} (x_j - x_i)(1 - x_i x_j) }
\det W^{2n}(\vectx ; \vecta)
\det W^{2n}(\vectx ; \vectb),
\label{P2}
\\
&
\Pf \left(
 \frac{ \det W^3(x_i, x_j, z ; a_i, a_j, c)
        \det W^2(x_i, x_j ; b_i, b_j) }
      { ( 1 - x_i x_j ) ( x_j - x_i ) }
\right)_{1 \le i, j \le 2n}
\notag
\\
&\qquad=
\frac{ (1+c)^{n-1} }
     { \prod_{1 \le i < j \le 2n} ( 1 - x_i x_j ) ( x_j - x_i ) }
\det W^{2n+1}( \vectx, z ; \vecta, c )
\det W^{2n}( \vectx ; \vectb ).
\hss
\label{P3}
\end{align}
\end{theorem}

We note that the $(i,j)$ entries of the determinant in (\ref{D2}) and
 the Pfaffian in (\ref{P2}) can be written in the form
\begin{align*}
&
\frac{ \det W^2(x_i, y_j ; a_i, b_j) }
     { (1 - x_i y_j)(y_j - x_i) }
=
 \frac{ 1 - a_i b_j }
      { 1 - x_i y_j }
+
 \frac{ b_j - a_i }
      { y_j - x_i },
\\
&
\frac{ \det W^2(x_i, x_j ; a_i, a_j) \det W^2(x_i, x_j ; b_i, b_j) }
     { (1 - x_i x_j)(x_j - x_i) }
\\
&\quad=
(1 - x_i x_j)(x_j - x_i)
\left(
 \frac{ 1 - a_i a_j }{ 1 - x_i x_j }
+
 \frac{ a_j - a_i }{ x_j - x_i }
\right)
\left(
 \frac{ 1 - b_i b_j }{ 1 - x_i x_j }
+
 \frac{ b_j - b_i }{ b_j - x_i }
\right).
\end{align*}

\begin{remark}
As a generalization of (\ref{D1}), we can show that
\begin{multline*}
\det \left(
 \frac{ \det V^{p+1,q+1}( x_i, y_j, \vectz ; a_i, b_j, \vectc ) }
      { y_j - x_i }
 \right)_{1 \le i, j \le n}
\\
=
\frac{ (-1)^{n(n-1)/2} }
     { \prod_{i,j=1}^n (y_j - x_i) }
\det V^{p,q}(\vectz ; \vectc)^{n-1}
\det V^{n+p,n+q}(\vectx,\vecty,\vectz;\vecta,\vectb,\vectc).
\end{multline*}
The identity (\ref{D1}) is the special case where $p=q=0$.
Also, we can prove
\begin{multline*}
\det \left(
 \frac{\det W^{p+2}(x_i,y_j,\vectz ; a_i,b_j,\vectc)}
      {(1 - x_i y_j) (y_j - x_i)}
 \right)_{1 \le i, j \le n}
\\
=
\frac{ 1 }{ \prod_{i,j=1}^n (1 - x_i y_j)(y_j - x_i) }
 \det W^p(\vectz ; \vectc)^{n-1}
 \det W^{2n+p}(\vectx,\vecty,\vectz ; \vecta,\vectb,\vectc),
\end{multline*}
which generalizes (\ref{D2}) and (\ref{D3}).
We also have Pfaffian identities, which are
 generalizations of (\ref{P1}), (\ref{P2}), and (\ref{P3}).
See \cite{IOT} for generalized identities, proofs and applications.
\end{remark}

\begin{demo}{Proof of Theorem~3.3}
The identities (\ref{D1}) and (\ref{D2}) are proven in
 \cite[Theorem~4.2, 4.3]{O}.
Here we give a proof of (\ref{D3}).

Let $I$ and $J$ be subsets of $[n] = \{ 1, 2, \cdots, n \}$,
 and let $L(I,J)$ (resp. $R(I,J)$) denote the coefficient of
 $a^I b^J = \prod_{i \in I} a_i \prod_{j \in J} b_j$
 on the left (resp. right) hand side of (\ref{D3}).
If we define the automorphisms $\sigma_I^x$ and $\sigma_J^y$ by setting
$$
\sigma_I^x(x_i) =
\begin{cases}
 x_i^{-1} &\text{if $i \in I$,} \\
 x_i &\text{if $i \not\in I$,}
\end{cases}
\quad\quad
\sigma_J^y(y_i) =
\begin{cases}
 y_j^{-1} &\text{if $j \in J$,} \\
 y_j &\text{if $j \not\in J$,}
\end{cases}
$$
then it follows from the definition of determinant that
\begin{align*}
\sigma_I^x \sigma_J^y \left( L(I,J) \right)
 &=
\det \left(
 \frac{ \det W^3(x_i,y_j,z ; 0, 0, c) }
      { (1 - x_i y_j)(y_j - x_i) }
\right)_{1 \le i, j \le n},
\\
\sigma_I^x \sigma_J^y \left( R(I,J) \right)
 &=
\frac{ (1+c)^{n-1} }
     { \prod_{i,j=1}^n ( 1 - x_i y_j ) ( y_j - x_i ) }
\det W^{2n+1}(\vectx, \vecty, z ; \vectzero, \vectzero, c),
\end{align*}
where $\vectzero = (0, \cdots, 0)$.
Hence it is enough to show
\begin{multline}
\det \left(
 \frac{ \det W^3(x_i,y_j,z ; 0, 0, c) }
      { (1 - x_i y_j)(y_j - x_i) }
\right)_{1 \le i, j \le n}
\\
=
\frac{ (1+c)^{n-1} }
     { \prod_{i,j=1}^n ( 1 - x_i y_j ) ( y_j - x_i ) }
\det W^{2n+1}(\vectx, \vecty, z ; \vectzero, \vectzero, c).
\label{D3-1}
\end{multline}
Regard the both sides of (\ref{D3-1}) as polynomials in $c$ and denote by
 $f(c)$ and $g(c)$ the left and right hand side of (\ref{D3-1})
 respectively.
Since $f(c)$ and $g(c)$ have degree at most $n$ in $c$,
 it is enough to prove the following three claims :
\begin{roster}{\bf Claim~3}
\item[\bf Claim~1.]
$f(c)$ is divisible by $(1+c)^{n-1}$.
\item[\bf Claim~2.]
$f(0) = g(0)$.
\item[\bf Claim~3.]
The coefficient of $c^n$ in $f(c)$ is equal to that in $g(c)$.
\end{roster}

First we prove Claim~1.
Let $A$ and $B$ be the $n \times n$ matrices with $(i,j)$ entry
 $z^2 - 1$ and $(1-zx_i)(1-zy_j)/(1-x_iy_j)$ respectively.
Then we have
$$
\left(
 \frac{ \det W^3(x_i,y_j,z;0,0,c) }
      { (1 - x_i y_j)(y_j - x_i) }
\right)_{1 \le i, j \le n}
=
A + (c+1) B.
$$
Here we use the following lemma.
(This lemma easily follows from the definition of determinant,
 so we omit its proof.)

\begin{lemma}
For $n \times n$ matrices $X$ and $Y$, we have
$$
\det (X + Y)
 =
\sum_{H,K}
 (-1)^{\Sigma(H) + \Sigma(K)}
 \det A_{H,K} \det B_{H^c,K^c},
$$
where the sum is taken over all pairs of subsets $H$ and $K \subset [n]$
 with $\# H = \# K$.
And $X_{H,K}$ (resp. $Y_{H^c,K^c}$) denotes the submatrix of
 $X$ (resp. $Y$) obtained by choosing entries with row indices in $H$
 (resp. $H^c$) and column indices $K$ (resp. $K^c$), and
 $\Sigma(H) = \sum_{h \in H}$, $\Sigma(K) = \sum_{k \in K} k$.
\end{lemma}

Applying this lemma and using the fact that $\rank A = 1$, we see that
$$
f(c) =
(c+1)^n \det B
 +
(c+1)^{n-1}
 \sum_{h,k=1}^n (-1)^{h+k} (z^2-1) \det B_{[n]- \{ h \},[n]- \{ k \} }.
$$
Therefore $f(c)$ is divisible by $(c+1)^{n-1}$.

Next we prove Claim~2.
It follows from the definition of determinant that
\begin{align*}
f(0)
 &=
\det \left(
 \frac{(z - x_i)(z - y_j)}{1 - x_i y_j}
\right)_{1 \le i, j \le n},
\\
g(0)
 &=
\frac{ 1 }
     { \prod_{i,j=1}^n (1 - x_i y_j)(y_j - x_i) }
\det W^{2n+1}(\vectx, \vecty, z ; \vectzero, \vectzero, 0).
\end{align*}
By using the Cauchy's determinant identity (\ref{C1}) and
 the Vandermonde determinant, we see that
$$
f(0) = g(0)
 =
\frac{ \prod_{i=1}^n (z - x_i)(z - y_i)
       \prod_{1 \le i < j \le n} (x_j - x_i)(y_j - y_i) }
     { \prod_{i,j=1}^n (1 - x_i y_j) }.
$$
Claim~3 can be proven similarly.
\end{demo}

\begin{demo}{Proof of Theorem~3.4}
The identities (\ref{P1}) and (\ref{P2}) are verified in
 \cite[Theorem~4.7, 4.4]{O}.
Here we give a proof of (\ref{P3}).

Let $I$ and $J$ be subsets of $[n] = \{ 1, 2, \cdots, n \}$,
 and let $L(I,J)$ (resp. $R(I,J)$) denote the coefficient of
 $a^I b^J = \prod_{i \in I} a_i \prod_{j \in J} b_j$
 on the left (resp. right) hand side of (\ref{P3}).
Let $\sigma_I$ be the automorphism defined by
$$
\sigma_I(x_i) =
\begin{cases}
 x_i^{-1} &\text{if $i \in I$,} \\
 x_i &\text{if $i \not\in I$.}
\end{cases}
$$
Then, by the same argument as in the proof of \cite[Theorem~4.4]{O},
 we can compute $\sigma_I(L(I,J))$ and $\sigma_I(R(I,J))$.
We put $K = (I \cap J^c) \cup (I^c \cap J)$ and define $Z(K)$
 to be the $2n \times 2n$ skew-symmetric matrix with $(i,j)$ entry
$$
Z(K)_{i,j} = \begin{cases}
-\dfrac{1}{1-x_ix_j} \det W^3(x_i,x_j,z;0,0,c)
 &\text{if $i \in K$ and $j \in K$,} \\
-\dfrac{1}{x_j-x_i} \det W^3(x_i,x_j,z;0,0,c)
 &\text{if $i \in K$ and $j \in K^c$,} \\
\dfrac{1}{x_j-x_i} \det W^3(x_i,x_j,z;0,0,c)
 &\text{if $i \in K^c$ and $j \in K$,} \\
\dfrac{1}{1-x_ix_j} \det W^3(x_i,x_j,z;0,0,c)
 &\text{if $i \in K^c$ and $j \in K^c$.}
\end{cases}
$$
Then we have
$$
\sigma_I( L(I,J) )
 = \prod_{i \in I} x_i^{-1} \cdot \Pf Z(K),
$$
and
$$
\sigma_I (R(I,J))
=
(-1)^{s(K)} \prod_{i \in I} x_i^{-1}
\cdot
\frac{ (c+1)^{n-1} }
     { \prod_{(k,l)} (1 - x_k x_l)
       \prod_{(k',l')} (x_{l'} - x_{k'}) }
\det W^{2n+1}( \vectx, z ; \vectzero, c ),
$$
where
$$
s(K) = \# \{ (i,j) \in K \times [2n] : i<j \},
$$
and the products are taken over all pairs $k < l$ and $k' < l'$
 such that
\begin{equation}
(k,l) \in (K \times K) \cup (K^c \times K^c),
\quad
(k',l') \in (K \times K^c) \cup (K^c \times K).
\label{range}
\end{equation}
Now the proof is reduced to showing the following identity :
\begin{equation}
\Pf Z(K)
=
(-1)^{s(K)}
\frac{ (c+1)^{n-1} }
     { \prod_{(k,l)} (1 - x_k x_l)
       \prod_{(k',l')} (x_{l'} - x_{k'}) }
\det W^{2n}(\vectx,z;0,c),
\label{Z}
\end{equation}
where the products are taken over all pairs $k < l$ and $k' < l'$
 satisfying (\ref{range}).

Regard the both sides of (\ref{Z}) as polynomials in $c$ and
 denote by $f(c)$ and $g(c)$ the left and right hand side respectively.
Since $f(c)$ and $g(c)$ have degree $n$, it is enough to show the following
 three claims :
\begin{roster}{\bf Claim~3}
\item[\bf Claim~1.]
$f(c)$ is divisible by $(c+1)^{n-1}$.
\item[\bf Claim~2.]
$f(0) = g(0)$.
\item[\bf Claim~3.]
The coefficient of $c^n$ in $f(c)$ is equal to that of $g(c)$.
\end{roster}

First we prove Claim~1.
Let $A(K)$ and $B(K)$ be $2n \times 2n$ skew-symmetric matrices
 with $(i,j)$ entry
\begin{gather*}
A(K)_{ij}
= \begin{cases}
 -(z^2-1) (x_j - x_i)   &\text{if $i \in K$ and $j \in K$,} \\
 -(z^2-1) (1 - x_i x_j) &\text{if $i \in K$ and $j \not\in K$,} \\
 (z^2-1) (1 - x_i x_j)  &\text{if $i \not\in K$ and $j \in K$,} \\
 (z^2-1) (x_j - x_i)    &\text{if $i \not\in K$ and $j \not\in K$,}
\end{cases}
\\
B(K)_{ij}
= \begin{cases}
 - (1 - z x_i)(1 - z x_j) \dfrac{x_j - x_i}{1 - x_i x_j}
  &\text{if $i \in K$ and $j \in K$,} \\
 - (1 - z x_i)(1 - z x_j)
  &\text{if $i \in K$ and $j \not\in K$,} \\
 (1 - z x_i)(1 - z x_j)
  &\text{if $i \not\in K$ and $j \in K$,} \\
 (1 - z x_i)(1 - z x_j) \dfrac{x_j - x_i}{1 - x_i x_j}
  &\text{if $i \not\in K$ and $j \not\in K$.} \\
\end{cases}
\end{gather*}
Then we have
$$
Z(K) = A(K) + (c+1) B(K).
$$
Here we use the following lemma.

\begin{lemma} (\cite[Lemma~4.2 (a)]{St}
If $X$ and $Y$ are $2n \times 2n$ skew-symmetric matrices, then we have
$$
\Pf (X + Y)
 = \sum_H (-1)^{\Sigma(H) - \# H/2} \Pf (X_H) \Pf (Y_{H^c}),
$$
where $H$ runs over all subsets $H \subset [2n]$ with $\# H$ even,
 and $X_H$ (resp. $Y_{H^c}$) denotes the skew-symmetric submatrix
 obtained from $X$
 (resp. $Y$) by picking the entries with row-indices and column-indices
 in $H$ (resp. $H^c$).
\end{lemma}

If $H$ is a subset with $\# H$ even, then we have
$$
\Pf \left( x_j - x_i \right)_{i,j \in H} = 0
\quad\text{if}\quad
\# H \ge 4,
$$
and, by applying the automorphism $\sigma_{H \cap K}$, we see that
$$
\Pf ( A(K)_H ) = 0.
\quad\text{if}\quad
\# H \ge 4.
$$
Hence, by using the above lemma, we have
\begin{align*}
\Pf Z(K)
 &=
 (c+1)^n \Pf B(K)
\\
&\quad +
 \sum_{1 \le k < l \le 2n}
  (-1)^{2n(2n+1)/2 - k-l-1} (c+1)^{n-1}
   \Pf ( A(K)_{k,l} ) \Pf ( B(K)_{[2n]-\{k,l\}} )
\end{align*}
Therefore $\Pf(Z(K))$ is divisible by $(c+1)^{n-1}$.

To prove Claims~2 and 3, we introduce the $2n \times 2n$
 skew-symmetric matrix $Y(K)$ by putting
$$
Y(K)_{ij} = \begin{cases}
 - \dfrac{x_j - x_i}{1 - x_i x_j} &\text{if $i \in K$ and $j \in K$,} \\
 - 1 &\text{if $i \in K$ and $j \not\in K$,} \\
 1 &\text{if $i \not\in K$ and $j \in K$,} \\
 \dfrac{x_j - x_i}{1 - x_i x_j} &\text{if $i \not\in K$ and $j \not\in K$.} \\
\end{cases}
$$
Then we have (see \cite[Lemma~4.5]{O})
\begin{equation}
\Pf Y(K)
 =
 (-1)^{s(K)}
 \prod_{\substack{k < l \\ k,l \in K}} \frac{x_l - x_k}{1 - x_k x_l}
 \prod_{\substack{k < l \\ k,l \not\in K}} \frac{x_l - x_k}{1 - x_k x_l}.
\label{Y}
\end{equation}
The constant term and the leading coefficient of $f(c) = \Pf Z(K)$
 are given by
\begin{align*}
[c^0] \Pf Z(K)
 &=
\Pf \left( (z - x_i)(z - x_j) Y(K)_{ij} \right)_{1 \le i, j \le 2n},
\\
[c^n] \Pf Z(K)
 &=
\Pf \left( (1 - zx_i)(1 - zx_j) Y(K)_{ij} \right)_{1 \le i, j \le 2n}.
\end{align*}
On the other hand, the constant term and the leading coefficient
 of $\det W^{2n}(\vectx, z ; \vectzero, c)$ are 
\begin{align*}
[c^0] \det W^{2n}(\vectx, z ; \vectzero, c)
 &=
\det W^{2n}(\vectx, z ; \vectzero, 0),
\\
[c^1] \det W^{2n}(\vectx, z ; \vectzero, c)
 &=
z^{2n} \det W^{2n}(\vectx, z^{-1} ; \vectzero, 0).
\end{align*}
Hence, by using (\ref{Y}) and the Vandermonde determinant, we see that
\begin{gather*}
[c^0] f(c) = [c^0] g(c)
 =
(-1)^{s(K)} \prod_{i=1}^{2n} (z - x_i)
 \prod_{\substack{k < l \\ k,l \in K}} \frac{x_l - x_k}{1 - x_k x_l}
 \prod_{\substack{k < l \\ k,l \in K^c}} \frac{x_l - x_k}{1 - x_k x_l},
\\
[c^n] f(c) = [c^n] g(c)
 =
(-1)^{s(K)} \prod_{i=1}^{2n} (1 - z x_i)
 \prod_{\substack{k < l \\ k,l \in K}} \frac{x_l - x_k}{1 - x_k x_l}
 \prod_{\substack{k < l \\ k,l \in K^c}} \frac{x_l - x_k}{1 - x_k x_l}.
\end{gather*}
This completes the proof of Theorem~3.4.
\end{demo}

The round determinants appearing the $2$- and $3$-enumerations can
 be evaluated by applying Theorem~3.3, but
 we need the substitution $\vecty = \vectx$ in the resulting determinants
 to obtain simple expressions,
 except for the $2$-enumeration of ASMs and VSASMs.
The following lemma will be used in this second step.
The proof is done by elementary transformations and left to
 the reader.

\begin{lemma}
Let $\alpha = (\alpha_1, \cdots, \alpha_{2n})$ be a sequence of
 half-integers,
and let $\vectx = (x_1, \cdots, x_n)$ and $\vecty =(y_1, \cdots, y_n)$
 be two vectors of $n$ variables.
\begin{roster}{(3)}
\item[(1)]
Let $V'(\alpha;\vectx,\vecty)$ be the $2n \times 2n$ matrix
 with $(i,j)$ entry
$$
\begin{cases}
 x_i^{\alpha_j}
   &\text{if $1 \le i \le n$,} \\
 (-1)^{j-1} y_{i-n}^{\alpha_j}
   &\text{if $n+1 \le i \le 2n$ and $1 \le j \le n$,} \\
 (-1)^{j-n} y_{i-n}^{\alpha_j}
   &\text{if $n+1 \le i \le 2n$ and $n+1 \le j \le 2n$.}
\end{cases}
$$
Then we have
$$
\det V'(\alpha ; \vectx, \vectx)
 =
(-1)^{n(n+1)/2} 2^n
 \det V(\beta ; \vectx) \det V(\beta' ; \vectx),
$$
where
$$
\beta = (\alpha_1, \alpha_{n+2}, \alpha_3, \alpha_{n+4}, \cdots),
\quad
\beta' = (\alpha_{n+1}, \alpha_2, \alpha_{n+3}, \alpha_4, \cdots).
$$
\item[(2)]
Let $W'_\pm(\alpha ; \vectx, \vecty)$ be the $2n \times 2n$ matrix
 with $(i,j)$ entry
$$
\begin{cases}
x_i^{\alpha_j} \pm x_i^{-\alpha_j}
 &\text{if $1 \le i \le n$,} \\
(-1)^{j-1} (y_{i-n}^{\alpha_j} \pm y_{i-n}^{-\alpha_j})
 &\text{if $n+1 \le i \le 2n$,}
\end{cases}
$$
Then we have
$$
\det W'_{\pm}(\alpha ; \vectx, \vectx)
 =
(-1)^{n(n+1)/2} 2^n
 \det W^{\pm}(\gamma ; \vectx) \det W^{\pm}(\gamma' ; \vectx),
$$
where
$$
\gamma = (\alpha_1, \alpha_3, \alpha_5, \cdots, \alpha_{2n-1}),
\quad
\gamma' = (\alpha_2, \alpha_4, \alpha_6, \cdots, \alpha_{2n}).
$$
\item[(3)]
Let $U(\alpha ; \vectx, \vecty)$ be the $2n \times 2n$ matrix with
 $(i,j)$ entry
$$
\begin{cases}
x_i^{\alpha_j} + x_i^{-\alpha_j} &\text{if $1 \le i \le n$,} \\
(-1)^{j-1} \left( y_{i-n}^{\alpha_j} - y_{i-n}^{-\alpha_j} \right)
 &\text{if $n+1 \le i \le 2n$.}
\end{cases}
$$
Then we have
$$
\det U(\alpha ; \vectx, \vectx)
= 2^n \det V(\tilde{\alpha} ; \vectx, \vectx^{-1}),
$$
where
$$
\tilde{\alpha}
 = (-\alpha_1, \alpha_2, -\alpha_3, \alpha_4, \cdots,
 -\alpha_{2n-1}, \alpha_{2n}).
$$
\end{roster}
\end{lemma}

\section{
Proof
}

In this section, we prove Theorem~2.3 and 2.4 stated in Section~2.
Since the arguments are the same, we illustrate how to compute
 the partition functions $A_{\VH}^{(2)}(4n+1 ; \vectx, \vecty ; \zeta_6)$
 and $A_{\VH}^{(2)}(4n+3 ; \vectx, \vecty ; \zeta_6)$,
 which correspond to the $1$-enumerations of VHSASMs.
(For other cases, see the end of this section and the tables there.)

First we consider the case of VHSASMs of order $4n+1$ and compute
 the partition function $A_{\VH}^{(2)}(4n+1 ; \vectx, \vecty ; \zeta_6)$.
A simple computation shows
\begin{multline*}
M_{\UU}(n ; \vectx, \vecty ; a, a, a)_{i,j}
\\
 =
\sigma(a) (1 - x_i^2) (1 - y_j^2)
\left(
 \frac{ x_i^2 + y_j^2 }
      { x_i^4 + y_j^4 - \dfrac{a^4+1}{a^2} x_i^2 y_j^2 }
+
 \frac{ 1 + x_i^2 y_j^2 }
      { x_i^4 y_j^4 + 1 - \dfrac{a^4+1}{a^2} x_i^2 y_j^2 }
\right).
\end{multline*}
If $a = \zeta_6$, then we have
$$
M_{\UU}(n ; \vectx, \vecty ; \zeta_6, \zeta_6, \zeta_6)
=
\sigma(\zeta_6)^n \prod_{i=1}^n (1 - x_i^2) (1 - y_i^2)
\det \left(
 \frac{ y_j^4 - x_i^4 }
      { y_j^6 - x_i^6 }
+
 \frac{ 1 - x_i^4 y_j^4 }
      { 1 - x_i^6 y_j^6 }
\right)_{1 \le i, j \le n}.
$$
Now, by applying the identity (\ref{D2}) in Theorem~3.3 with
$$
x_i \to x_i^6,
\quad
y_i \to y_i^6,
\quad
a_i \to x_i^4,
\quad
b_i \to y_i^4,
$$
we have
$$
\det \left(
 \frac{ y^4 - x^4 }
      { y^6 - x^6 }
+
 \frac{ 1 - x^4 y^4 }
      { 1 - x^6 y^6 }
\right)_{1 \le i, j \le n}
=
\frac{ 1 }
     { \prod_{i,j=1}^n (y_j^6 - x_i^6)(1 - x_i^6 y_j^6) }
\det W^{2n}(\vectx^6, \vecty^6 ; \vectx^4, \vecty^4).
$$
By applying elementary transformations and by using the definition of
 orthogonal characters (\ref{O2}), we have
\begin{align*}
&
\det W^{2n}(\vectx^6, \vecty^6 ; \vectx^4, \vecty^4)
\\
&\quad=
\prod_{i=1}^n x_i^{2n+1} y_i^{2n+1}
\prod_{1 \le i < j \le n}
 (x_j^2 - x_i^2)(1 - x_i^2 x_j^2)
 (y_j^2 - y_i^2)(1 - y_i^2 y_j^2)
\prod_{i,j=1}^n
 (y_j^2 - x_i^2)(1 - x_i^2 y_j^2)
\\
&\quad\quad\times
\tilde{\Orth}_{4n}( \delta(n+1/2,n-1/2) ; \vectx^2,\vecty^2).
\end{align*}
Then, after some computation, we obtain
$$
A_{\VH}^{(2)}(4n+1 ; \vectx, \vecty ; \zeta_6)
=
3^{-n^2}
\frac{ 1 }
     { \prod_{i=1}^n (x_i + x_i^{-1}) (y_i + y_i^{-1}) }
\tilde{\Orth}_{4n}( \delta(n+1/2,n-1/2) ; \vectx^2,\vecty^2).
$$

Next we consider the cases of VHSASMs of order $4n+3$.
A simple computation shows that
\begin{align*}
&
M_{\UU}(n ; \vectx, \vecty ; a, a^{-1}, a^{-1})_{i,j}
\\
&\quad=
\sigma(a)
 \left(
 \frac{ x_i^4 y_j^2 + y_j^2 + x_i^2 y_j^4 + x_i^2
        - \dfrac{a^4+a^2+1}{a^2} (x_i^4 + y_j^4)
        + 2 x_i^2 y_j^2
      }
      { x_i^4 + y_j^4 - \dfrac{a^4+1}{a^2} x_i^2 y_j^2
      }
 \right.
\\
&\quad\quad\quad\quad\quad
 \left.
 -
 \frac{ x_i^4 y_j^2 + y_j^2 + x_i^2 y_j^4 + x_i^2
        - \dfrac{a^4+a^2+1}{a^2} (x_i^4 y_j^4 + 1)
        + 2 x_i^2 y_j^2
      }
      { x_i^4 y_j^4 + 1 - \dfrac{a^4+1}{a^2} x_i^2 y_j^2
      }
 \right).
\end{align*}
If $a = \zeta_6$, then $(a^4+a^2+1)/a^2 = 0$ and we have
\begin{multline}
\det M_{\UU}(n ; \vectx, \vecty ; \zeta_6, \zeta_6^{-1}, \zeta_6^{-1})
\\
=
\sigma(\zeta_6)^n
\prod_{i=1}^n (1 - x_i^4)(1 - y_i^4)
\det \left(
\frac{ x_i^4 y_j^2 + y_j^2 + x_i^2 y_j^4 + x_i^2 + 2 x_i^2 y_j^2 }
     { (x_i^4 + y_j^4 + x_i^2 y_j^2) (x_i^4 y_j^4 + 1 + x_i^2 y_j^2) }
\right)_{1 \le i, j \le n}.
\label{muu}
\end{multline}
Here we note that the numerator of the $(i,j)$ entry is equal to
$$
x_i^4 y_j^2 + x_i^2 y_j^4 + 2 x_i^2 y_j^2 + x_i^2 + y_j^2
 =
x_i^2 y_j^2 \cdot \Symp_6( (1,0,0) ; x_i^2, y_j^2, 1 ).
$$
Instead of evaluating directly the determinant on the right hand side
 of (\ref{muu}), we consider the determinant
$$
\det \left(
 \frac{ x_i^2 y_j^2 z^2 \Symp_6( (1,0,0) ; x_i^2, y_j^2, z^2 ) }
      { (x_i^4 + y_j^4 + x_i^2 y_j^2) (x_i^4 y_j^4 + 1 + x_i^2 y_j^2) }
\right)_{1 \le i, j \le n}.
$$
Comparing the definitions of the symplectic character
 $\Symp_6( (1,0,0) )$ and the matrix $W^3$, we have
$$
\det W^3(x^6, y^6, z^6 ; -x^4, -y^4, -z^4)
=
x^8 y^8 z^8 \cdot \Symp_6( (1,0,0) ; x^2,y^2,z^2)
\det W^-(1,2,3;x^2,y^2,z^2).
$$
Hence we see that
\begin{align*}
&
\det \left(
 \frac{ x_i^2 y_j^2 z^2 \cdot \Symp_6( (1,0,0) ; x_i^2, y_j^2, z^2 ) }
      { (x_i^4 + y_j^4 + x_i^2 y_j^2) (x_i^4 y_j^4 + 1 + x_i^2 y_j^2) }
\right)_{1 \le i, j \le n}
\\
&\quad=
\frac{ 1 }
     { (1-z^4)^n
       \prod_{i=1}^n
         (1 - x_i^4)(1 - y_i^4)
         (1 - x_i^2 z^2)(1 - y_i^2 z^2)
         (z^2 - x_i^2)(z^2 - y_i^2)
     }
\\
&\quad\quad\times
\det \left(
 \frac{ \det W^3(x_i^6, y_j^6, z_j^6 ; -x_i^4, -y_j^4, -z_j^4) }
      { ( 1 - x_i^6 y_j^6 ) ( y_j^6 - x_i^6 ) }
\right)_{1 \le i, j \le n}.
\end{align*}
Now we can apply the identity (\ref{D3}) in Theorem~3.3 with
$$
x_i \to x_i^6,
\quad
y_j \to y_j^6,
\quad
z \to z^6,
\quad
a_i \to - x_i^4,
\quad
b_j \to - y_j^4,
\quad
c \to - z^4.
$$
By applying elementary transformations and by using the definition (\ref{Sp})
 of symplectic characters, we have
\begin{align*}
&
\det W^{2n+1}(\vectx^6,\vecty^6,z^6 ; -\vectx^4, -\vecty^4,-z^4)
\\
&\quad=
z^{2n} (1-z^4)
\prod_{i=1}^n x_i^{2n} y_i^{2n}
\prod_{i=1}^n
 (1 - x_i^4)(1 - y_i^4)
 (z^2 - x_i^2)(z^2 - y_i^2) (1 - x_i^2 z^2)(1 - y_i^2 z^2)
\\
&\quad\quad\times
\prod_{1 \le i < j \le n}
 (x_j^2 - x_i^2) (1 - x_i^2 x_j^2)
 (y_j^2 - y_i^2) (1 - y_i^2 y_j^2)
\prod_{i,j=1}^n (y_j^2 - x_i^2) (1 - x_i^2 y_j^2)
\\
&\quad\quad\times
\Symp_{4n+2}( \delta(n,n-1) ; \vectx^2, \vecty^2, \vectz^2).
\end{align*}
After canceling the common factors, we can substitute $z = 1$ and
 obtain
\begin{align*}
&
\det \left(
\frac{ x_i^4 y_j^2 + x_i^2 y_j^4 + 2 x_i^2 y_j^2 + x_i^2 + y_j^2
     }
     { (1 + x_i^2 y_j^2 + x_i^4 y_j^4) (x_i^4 + x_i^2 y_j^2 + y_j^4) }
\right)_{1 \le i, j \le n}
\\
&\quad=
\frac{ 1 }
     { \prod_{i,j=1}^n
         ( 1 + x_i^2 y_i^2 + x_i^4 y_j^4 ) ( x_i^4 + x_i^2 y_j^2 + y_j^4 )
     }
\\
&\quad\quad\times
\prod_{i=1}^n x_i^{2n} y_i^{2n}
\prod_{1 \le i < j \le n}
 (x_j^2 - x_i^2)(1 - x_i^2 x_j^2)
 (y_j^2 - y_i^2)(1 - y_i^2 y_j^2)
\cdot
\Symp_{4n+2}( \delta(n,n-1) ; \vectx^2, \vecty^2, 1).
\end{align*}
Then, by some computation, we have
$$
A_{\VH}^{(2)}(4n+3 ; \vectx, \vecty ; \zeta_6)
=
3^{-n^2}
\Symp_{4n+2}( \delta(n,n-1) ; \vectx^2, \vecty^2, 1).
$$
This completes the proof of Theorem~2.3 (4) and (5) at $a = \zeta_6$.

The determinants and Pfaffians in the $0$-enumerations (or
 in the case of $a = \zeta_4$) are evaluated by families of
 the Cauchy's identity and the Borchardt's identity.
The determinant/Pfaffian identities used in the evaluation are listed
 in Table~1.

\begin{table}[htbp]
\caption{$0$-enumeration}
\begin{center}
\begin{tabular}{c|c|c}
Partition functions & Identities & specialization \\
\hline
$A(n ; \vectx, \vecty ; \zeta_4)$ &
 (\ref{B1}) & $x_i \to x_i^2$, $y_i \to y_i^2$ \\
\hline
$A_{\HT}^{(2)}(n ; \vectx, \vecty ; \zeta_4)$ &
 (\ref{C1}) & $x_i \to x_i^2$, $y_i \to y_i^2$ \\
\hline
$A_{\V}(2n+1 ; \vectx, \vecty ; \zeta_4)$ &
 (\ref{B1}) & $x_i \to x_i^2 + x_i^{-2}$, $y_i \to y_i^{-2}$ \\
\hline
$A_{\VHS}^{(2)}(4n+1 ; \vectx, \vecty ; \zeta_4)$ &
 (\ref{C1}) & $x_i \to x_i^2 + x_i^{-2}$, $y_i \to y_i^{-2}$ \\
\hline
$A_{\VHS}^{(2)}(4n+3 ; \vectx, \vecty ; \zeta_4)$ &
 (\ref{C1}) & $x_i \to x_i^2 + x_i^{-2}$, $y_i \to y_i^{-2}$ \\
\hline
$A_{\UU}^{(2)}(4n ; \vectx, \vecty ; \zeta_4,\zeta_4,\zeta_4)$ &
 (\ref{C1}) & $x_i \to x_i^2 + x_i^{-2}$, $y_i \to y_i^{-2}$ \\
\hline
$A_{\VHP}^{(2)}(4n+2 ; \vectx, \vecty ; \zeta_4)$ &
 (\ref{C1}) & $x_i \to x_i^2 + x_i^{-2}$, $y_i \to y_i^{-2}$ \\
\hline
$A_{\QT}^{(1)}(4n ; \vectx ; \zeta_4)$ &
 (\ref{I1}) & $x_i \to x_i^2$ \\
\hline
$A_{\QT}^{(2)}(4n ; \vectx ; \zeta_4)$ &
 (\ref{S1}) & $x_i \to x_i^2$ \\
\hline
$A_{\OD}(2n ; \vectx ; \zeta_4)$ &
 (\ref{I2}) & $x_i \to \sqrt{-1} x_i^2$ \\
\hline
$A_{\OO}^{(2)}(4n ; \vectx ; \zeta_4)$ &
 (\ref{S2}) & $x_i \to \sqrt{-1} x_i^2$ \\
\hline
$A_{\UO}^{(1)}(8n ; \vectx ; \zeta_4, \zeta_4)$ &
 (\ref{I1}) & $x_i \to x_i^2 + x_i^{-2}$ \\
\hline
$A_{\UO}^{(2)}(8n ; \vectx ; \zeta_4, \zeta_4)$ &
 (\ref{S1}) & $x_i \to x_i^2 + x_i^{-2}$ \\
\hline
$A_{\VO}^{(2)}(8n+1 ; \vectx ; \zeta_4)$ &
 (\ref{S1}) & $x_i \to x_i^2 + x_i^{-2}$ \\
\hline
$A_{\VO}^{(2)}(8n+3 ; \vectx ; \zeta_4)$ &
 (\ref{S1}) & $x_i \to x_i^2 + x_i^{-2}$ \\
\end{tabular}
\end{center}
\end{table}

When we compute the partition functions in the $1$-enumerations,
 we apply the identities in Theorem~3.3 and 3.4 to evaluate the determinants
 or Pfaffians, and compare the resulting determinants with definitions
 (\ref{GL})--(\ref{O2}) of classical group characters.
The variations of the arguments are summarized in Table~2.

\begin{table}[htbp]
\caption{$1$-enumeration}
\begin{center}
\begin{tabular}{c|c|c}
Partition functions & Identities & specialization \\
\hline
$A(n ; \vectx, \vecty ; \zeta_6)$ &
 (\ref{D1}) &
 $x_i \to x_i^6$, $a_i \to x_i^2$, $y_i \to y_i^6$, $b_i \to y_i^2$
\\
\hline
$A_{\HT}^{(2)}(n ; \vectx, \vecty ; \zeta_6)$ &
 (\ref{D1}) &
$x_i \to x_i^6$, $a_i \to x_i^4$, $y_i \to y_i^6$, $b_i \to y_i^4$
\\
\hline
$A_{\V}(2n+1 ; \vectx, \vecty ; \zeta_6)$ &
 (\ref{D2}) &
$x_i \to x_i^6$, $a_i \to -x_i^2$, $y_i \to y_i^6$, $b_i \to -y_i^2$
\\
\hline
$A_{\VHS}^{(2)}(4n+1 ; \vectx, \vecty ; \zeta_6)$ &
 (\ref{D2}) &
$x_i \to x_i^6$, $a_i \to x_i^4$, $y_i \to y_i^6$, $b_i \to y_i^4$
\\
\hline
$A_{\VHS}^{(2)}(4n+3 ; \vectx, \vecty ; \zeta_6)$ &
 (\ref{D3}) &
$x_i \to x_i^6$, $a_i \to -x_i^4$, $y_i \to y_i^6$, $b_i \to -y_i^4$,
$z \to z^6$, $c \to -z^4$
\\
\hline
$A_{\VHP}^{(2)}(4n+2 ; \vectx, \vecty ; \zeta_6)$ &
 (\ref{D2}) &
$x_i \to x_i^6$, $a_i \to -x_i^2$, $y_i \to y_i^6$, $b_i \to -y_i^2$
\\
\hline
$A_{\QT}^{(1)}(4n ; \vectx ; \zeta_6)$ &
 (\ref{P1}) &
 $x_i \to x_i^6$, $a_i \to x_i^2$, $b_i \to x_i^2$
\\
\hline
$A_{\QT}^{(2)}(4n ; \vectx ; \zeta_6)$ &
 (\ref{P1}) &
 $x_i \to x_i^6$, $a_i \to x_i^2$, $b_i \to x_i^4$
\\
\hline
$A_{\OD}(2n ; \vectx ; \zeta_6)$ &
 (\ref{P2}) &
 $x_i \to x_i^6$, $a_i \to -x_i^2$, $b_i \to 0$
\\
\hline
$A_{\UO}^{(1)}(8n ; \vectx ; \zeta_6, \zeta_4)$ &
 (\ref{P2}) &
 $x_i \to x_i^6$, $a_i \to -x_i^2$, $b_i \to -x_i^2$
\\
\hline
$A_{\VO}^{(2)}(8n+1 ; \vectx ; \zeta_6)$ &
 (\ref{P2}) &
 $x_i \to x_i^6$, $a_i \to -x_i^2$, $b_i \to x_i^4$
\\
\hline
$A_{\VO}^{(2)}(8n+3 ; \vectx ; \zeta_6)$ &
 (\ref{P3}) &
 $x_i \to x_i^6$, $a_i \to -x_i^4$, $b_i \to -x_i^2$, $z \to z^6$, $c \to -z^4$
\\
\end{tabular}
\end{center}
\end{table}

For the $2$-enumerations, the determinant appearing in
 $A(n ; \vectx, \vecty ; \zeta_8)$ (resp.
 $A_{\V}(2n+1 ; \vectx, \vecty ; \zeta_8)$)
 can be evaluated by applying the Cauchy's determinant identity
 with $x_i \to x_i^2$ and $y_i \to y_i^2$ (resp.
 $x_i$ by $x_i^2 + x_i^{-2}$ and $y_i \to y_i^2 + y_i^{-2}$).
The other determinants in the determinant partition functions are
 computed by using the identities (\ref{D1}) and (\ref{D2}),
 and then by applying Lemma~3.8 to the resulting determinants,
 where we have to put $\vecty = \vectx$.
Also the round Pfaffian partition functions for $a = \zeta_8$
 are computed by using the identities (\ref{S1}), (\ref{P1}) and (\ref{P2}).
Similarly, we can deal with the round $3$-enumerations.
See Table~3 and 4 for the details.

\begin{table}[htbp]
\caption{$2$-enumeration}
\begin{center}
\begin{tabular}{c|c|c|c}
Partition functions & Identities & specialization & Lemma~3.8 \\
\hline
$A(n ; \vectx, \vecty ; \zeta_8)$ &
 (\ref{C1}) & $x_i \to x_i^2$, $y_i \to y_i^2$ & \\
\hline
$A_{\HT}^{(2)}(n ; \vectx, \vecty ; \zeta_8)$ &
 (\ref{D1}) &
 $x_i \to x_i^4$, $a_i \to x_i^2$, $y_i \to -y_i^4$, $b_i \to -y_i^2$ &
 (1) \\
\hline
$A_{\V}(2n+1 ; \vectx, \vecty ; \zeta_8)$ &
 (\ref{C1}) & $x_i \to x_i^2 + x_i^{-2}$, $y_i \to y_i^2 + y_i^{-2}$ & \\
\hline
$A_{\VHS}^{(2)}(4n+1 ; \vectx, \vecty ; \zeta_8)$ &
 (\ref{D2}) &
 $x_i \to x_i^4$, $a_i \to x_i^2$, $y_i \to -y_i^4$, $b_i \to -y_i^2$ &
 (2) \\
\hline
$A_{\VHS}^{(2)}(4n+3 ; \vectx, \vecty ; \zeta_8)$ &
 (\ref{D2}) &
 $x_i \to x_i^4$, $a_i \to -x_i^2$, $y_i \to -y_i^4$, $b_i \to y_i^2$ &
 (2) \\
\hline
$A_{\UU}^{(2)}(4n ; \vectx, \vecty ; \zeta_8,\zeta_4,\zeta_4)$ &
 (\ref{D2}) &
 $x_i \to x_i^4$, $a_i \to -x_i^2$, $y_i \to -y_i^4$, $b_i \to y_i^2$ &
 (2) \\
\hline
$A_{\VHP}^{(2)}(4n+2 ; \vectx, \vecty ; \zeta_8)$ &
 (\ref{D2}) &
 $x_i \to x_i^4$, $a_i \to x_i^2$, $y_i \to -y_i^4$, $b_i \to y_i^2$ &
 (3) \\
\hline
$A_{\QT}^{(1)}(4n ; \vectx ; \zeta_8)$ &
 (\ref{P1}) &
 $x_i \to x_i^8$, $a_i \to x_i^2$, $b_i \to x_i^4$ & \\
\hline
$A_{\QT}^{(2)}(4n ; \vectx ; \zeta_8)$ &
 (\ref{S1}) & $x_i \to x_i^2$ & \\
\end{tabular}
\end{center}
\end{table}

\begin{table}[htbp]
\caption{$3$-enumeration}
\begin{center}
\begin{tabular}{c|c|c|c}
Partition functions & Identities & specialization & Lemma~3.8 \\
\hline
$A(n ; \vectx, \vecty ; \zeta_{12})$ &
 (\ref{D1}) &
 $x_i \to x_i^6$, $a_i \to x_i^2$, $y_i \to -y_i^6$, $b_i \to -y_i^2$ &
 (1) \\
\hline
$A_{\V}(2n+1 ; \vectx, \vecty ; \zeta_{12})$ &
 (\ref{D2}) &
 $x_i \to x_i^6$, $a_i \to -x_i^2$, $y_i \to -y_i^6$, $b_i \to y_i^2$ &
 (2) \\
\hline
$A_{\QT}^{(1)}(4n ; \vectx ; \zeta_{12})$ &
 (\ref{P1}) &
 $x_i \to x_i^{12}$, $a_i \to x_i^4$, $b_i \to x_i^6$ & \\
\end{tabular}
\end{center}
\end{table}

\section{
Discussion
}

In this paper, we settled a conjecture on the number of VHSASMs.
However, the enumeration problems of odd-order HTSASMs, odd-order QTSASMs,
 and odd-order DSASMs are still open.
In our point of view, the remaining conjectures (see \cite{R}) on HTSASMs and
 DSASMs are reformulated as follows:

\begin{conjecture}
\begin{roster}{(2)}
\item[(1)]
The number of $(2n+1) \times (2n+1)$ HTSASMs is given by
$$
\# \ASM^{\HTS}_{2n+1} =
3^{-n^2} \left( \dim \GL_{2n+1}(\delta(n,n-1)) \right)^2.
$$
\item[(2)]
The number of the $(2n+1) \times (2n+1)$ DSASMs is given by
$$
\# \ASM^{\DS}_{2n+1} =
3^{-n(n-1)/2} \dim \GL_{2n+1}(\delta(n,n-1)).
$$
\end{roster}
\end{conjecture}

Our results (Theorem~1.2) says that, for example, there is
 a bijection between the set of $n \times n$ ASMS
 and the set of all pairs $(T,M)$ of semistandard tableaux $T$
 of shape $\delta(n-1,n-1)$ with entries $1, 2, \cdots, 2n$, 
and triangular array $M$ of $1$s, $0$s and $-1$s
 of order $n$.
It would be interesting to find bijections proving the formulae in Theorem~1.2.
Also it is important to clarify the intrinsic reason why classical group
 characters appear in the enumeration of symmetry classes of ASMs.

\end{document}